\documentclass[preprint,12pt]{elsarticle}
\usepackage[left=0.75in, right=1in, top=1in, bottom=1in]{geometry}
\usepackage{graphicx}
\usepackage{amssymb}
\usepackage{amsthm}
\usepackage[tbtags]{amsmath}
\usepackage{mathtools}
\numberwithin{equation}{section}
\usepackage{lineno}
\usepackage[font=small,labelfont=bf]{caption}
\usepackage{subcaption}
\usepackage{romannum}
\usepackage{array}
\usepackage{booktabs}
\usepackage{float}
\usepackage{cleveref}
\usepackage{bm}
\usepackage{nicefrac}
\usepackage{color}
\usepackage{calrsfs}

\def\hf{\frac{1}{2}}

\newcommand{\dt}{\Delta t}
\newcommand{\dx}{\Delta x}
\newcommand{\dy}{\Delta y}
\newcommand{\kph}{{k+\frac{1}{2}}}
\newcommand{\kmh}{{k-\frac{1}{2}}}
\newcommand{\jph}{{j+\frac{1}{2}}}
\newcommand{\jmh}{{j-\frac{1}{2}}}
\newcommand{\cU}{\mathcal{U}}

\newcommand*\xbar[1]{%
  \hbox{%
    \vbox{%
      \hrule height 0.5pt 
      \kern0.4ex
      \hbox{%
        \kern-0.05em
        \ensuremath{#1}%
        \kern-0.05em
      }%
    }%
  }%
}

\makeatletter
\def\ps@pprintTitle{%
   \let\@oddhead\@empty
   \let\@evenhead\@empty
   \let\@oddfoot\@empty
   \let\@evenfoot\@oddfoot
}
\makeatother

\begin{document}
\begin{frontmatter}
\pagenumbering{arabic}

\title{An adaptive central-upwind scheme on quadtree grids for variable density shallow water equations}

\author[]{Mohammad A. Ghazizadeh\corref{cor1}}
\ead{sghaz023@uottawa.ca}

\author[]{Abdolmajid Mohammadian}
\ead{majid.mohammadian@uottawa.ca}
 
\cortext[cor1]{Corresponding author.}

\address{Department of Civil Engineering, University of Ottawa, Ottawa, ON K1N 6N5, Canada}

\begin{abstract}
Minimizing computational cost is one of the major challenges in the modelling and numerical analysis of hydrodynamics, and one of the ways to achieve this is by the use of quadtree grids.
In this paper, we present an adaptive scheme on quadtree 
grids for variable density shallow water equations. A scheme for the coupled system is developed based on the well-balanced positivity-preserving central-upwind scheme proposed in \cite{Ghazizadeh2019}. The scheme is capable of exactly preserving ``lake-at-rest'' steady states. A continuous piecewise bi-linear interpolation of the bottom topography function is used 
to achieve higher-order in space in order to preserve the positivity of water depth for the point values of each computational cell. Necessary conditions are checked to be able to preserve the positivity of water depth and density, and to ensure the achievement of a stable numerical scheme.
At each timestep, local gradients are examined to find new seeding points to locally refine/coarsen the computational grid.
\end{abstract}
\end{frontmatter}
\section{Introduction}
Quadtree grids are two-dimensional (2-D) semi-structured Cartesian grids that can be very accommodating for various problems in the field of computational hydrodynamics. One of the advantages of quadtree grids over structured and 
unstructured grids is grid coarsening/refining. The accuracy is increased/maintained while the grid refines/coarsens wherever it is needed 
and thus, the computational cost is reduced. There are a number of studies on how to generate quadtree grids; see, e.g., 
\cite{Aizawa2008,Borthwick2000,Ghazizadeh2019,Greaves1998,Pascal1998,Popinet2010,Samet1984a,Samet2006}.

The main goal of this paper is to develop an adaptive well-balanced positivity-preserving central-upwind 
scheme on quadtree grids for the coupled 
variable density shallow water equations (SWEs). The variable density SWEs in 2-D can be written in terms of 
conservative variables of $w$ (water surface), $hu$ and $hv$ (the unit discharges), and $h\rho$:
\begin{equation}
\begin{dcases}
w_t+(hu)_x+(hv)_y=0,\\
(hu)_t+\big( hu^2+\dfrac{g}{2\rho_\circ}h^2\rho\big)_x+(huv)_y=-\dfrac{g}{\rho_\circ} h\rho B_x,\\
(hv)_t+(huv)_x+\big(hv^2+\dfrac{g}{2\rho_\circ} h^2\rho\big)_y=-\dfrac{g}{\rho_\circ} h\rho B_y,\\
(h\rho)_t+(hu\rho)_x+(hv\rho)_y=0,
\end{dcases}
\label{eq:SWE}
\end{equation}
where $t$ is time, $g$ is the gravitational constant, $x$ and $y$ are the directions in the 2-D Cartesian coordinate system, $u(x,y,t)$ and
$v(x,y,t)$ are the water velocities in the $x$- and $y$-directions, respectively, $B(x,y)$ is the bottom topography,
$h(x,y,t)=w(x,y,t)-B(x,y)$ is the water depth, $\rho$ is the density, and $\rho_\circ$ is the reference density.

System \eqref{eq:SWE} admits ``lake-at-rest'' steady-state solutions,
\begin{equation}
  \rho\equiv{\rm Const},\quad w\equiv{\rm Const},\quad u=v\equiv0, \quad B\equiv{\rm Const},
  \label{1.2}
  \end{equation}
which can be obtained from \eqref{eq:SWE} \cite{Chertock2014}. The following quadtree scheme is 
capable of exactly preserving ``lake-at-rest'' steady states, which is called the {\em well-balanced} property. Another important attribute of the following method is its ability to preserve the non-negativity of $h$ 
and $\rho$, which is called the {\em positivity-preserving} property (see \cite{Kur_Acta} for a comprehensive review on these subjects).

A number of numerical schemes on quadtree grids for the SWEs have been introduced in recent years. For example, an adaptive well-balanced 
positivity-preserving central-upwind high-order scheme on quadtree grids was proposed in \cite{Ghazizadeh2019}. In addition, a well-balanced
scheme on quadtree-cut-cell grids was proposed in \cite{An2012}. This scheme is based on the hydrostatic reconstruction from \cite{Audusse2005}.
Furthermore, an adaptive second-order Roe-type scheme was proposed in \cite{Rogers2001}. An adaptive well-balanced Godunov-type scheme for the shallow water for the wet-dry over complex topography was introduced in \cite{Liang2009} and an adaptive quadtree Roe-type scheme for 2-D two-layer SWEs was presented in \cite{Lee2011}. Besides the
aforementioned numerical methods, several well-balanced positivity-preserving central-upwind schemes for the shallow water equations have been proposed 
in the past years; see, e.g., \cite{Audusse2004,Audusse2005,BMK,BM08,BCKN,BNL,Bryson2011,Chertock2008,Chertock2018,GPC,Kurganov2002,Kurganov2007,LAEK,Ric15,SMSK}, 
yet, to our knowledge, none of these methods has been extended to the coupled variable density SWEs over quadtree grids.

In \cite{Jiang2011} the coupled variable density SWEs were studied with a Godunov-type HLLC approximate Riemann 
solver. There are other studies that have been conducted on variable density SWEs and variable horizontal temperature SWEs (which 
have mathematically similar properties) with different numerical schemes; see, e.g. \cite{Chertock2014, Guo2011, Khorshid2017, Leighton2009}.

In this paper, we propose a central-upwind quadtree scheme which is based on the one from \cite{Ghazizadeh2019}.
Central-upwind schemes are finite-volume methods that are Godunov-type Riemann-problem-solver-free 
\cite{KLin,KNP,KPW,KTcl}. Central-upwind schemes have been referred to as ``black-box'' solvers
for general multidimensional systems of hyperbolic systems of conservation laws, and hav been extended to shallow water
models \cite{Kur_Acta}. The proposed scheme is the first well-balanced positivity-preserving central-upwind scheme for the variable density SWEs over quadtree grids.
This method is simple, efficient, and robust.

The paper is organized as follows. In \S\ref{S:2}, we briefly describe the quadtree grid generation terminology. In \S\ref{S:3}, we construct a
central-upwind quadtree scheme for the variable density SWEs with the mentioned features and test it on four different numerical examples in \S\ref{S:4}. 
Finally in \S\ref{S:5}, some concluding remarks are presented.

\section{Quadtree grids}\label{S:2}
In this section, we denote the terminology for 
how to generate quadtree grids (see \cite{Ghazizadeh2019, Borthwick2000,Greaves1998,Popinet2003}):

\vskip5pt
\noindent\textbf{Seeding points:} A set of points that helps to locally refine/coarsen the computational grid when needed. 

\vskip4pt
\noindent\textbf{Level of refinement:} Level of the quadtree, in which the size of the smallest cell is inversely
proportional to the maximum level of $m$.

\begin{figure}[h!]
\centerline{\includegraphics[height=0.30\linewidth]{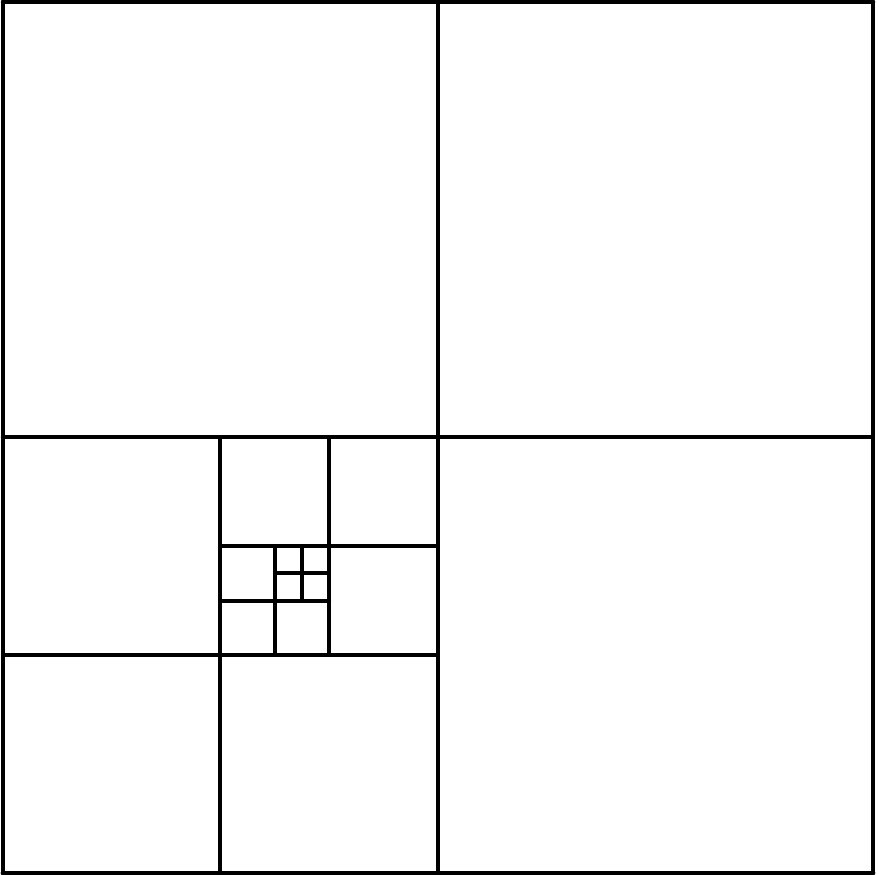}\hspace*{0.5cm}\includegraphics[height=0.30\linewidth]{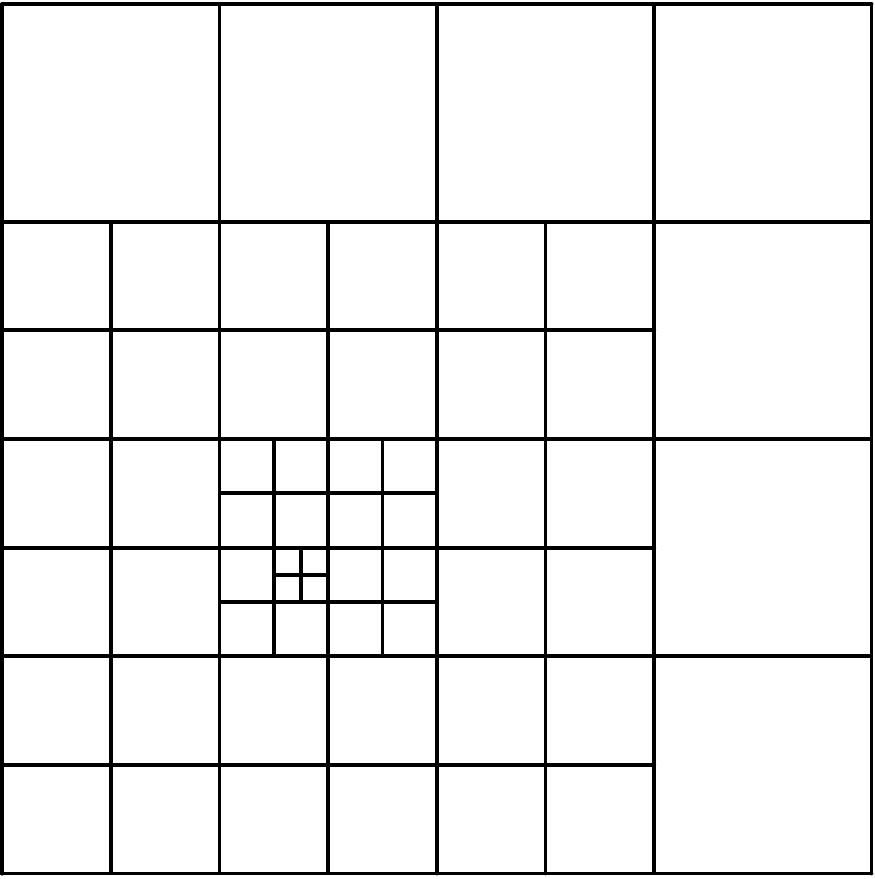}}
\caption{\sf Examples of non-regularised (left) and regularised (right) quadtree grids.\label{fig:1}}
\end{figure}

\vskip4pt
\noindent\textbf{Regularised quadtree grid:} In a regularised grid, no cell can have both an adjacent
neighboring cell and a diagonally neighboring cell with a refinement level difference greater than one (Figure \ref{fig:1}). 
The proposed scheme is based on regularised quadtree grids to prevent complicated formulations and improve stability.

\section{Adaptive semi-discrete central-upwind scheme}\label{S:3}
We write system \eqref{eq:SWE} in the following vector form:
\begin{equation}
\bm{U}_t+\bm{F}(\bm{U},B)_x+\bm{G}(\bm{U},B)_y=\bm{S}(\bm{U},B),
\label{3.1}          
\end{equation}
where
\begin{equation*}
\bm{U}:=(w,hu,hv, h\rho)^\top,
\end{equation*}
and the fluxes and source term are:
\begin{align}
&\bm{F}(\bm{U},B)=\left(hu,\dfrac{(hu)^2}{w-B}+\dfrac{g}{2\rho_\circ}\rho\,(w-B)^2,\dfrac{(hu)(hv)}{w-B}, hu\rho \right)^\top,\label{3.2}\\
&\bm{G}(\bm{U},B)=\left(hv,\dfrac{(hu)(hv)}{w-B},\dfrac{(hv)^2}{w-B}+\dfrac{g}{2\rho_\circ}\rho\,(w-B)^2, hv\rho \right)^\top,\label{3.3}\\ 
&\bm{S}(\bm{U},B)=\left(0,-\dfrac{g}{\rho_\circ} \rho\,(w-B)B_x,-\dfrac{g}{\rho_\circ} \rho\,(w-B)B_y, 0 \right)^\top.\label{3.4}
\end{align}

In the following, an adaptive well-balanced semi-discrete central-upwind scheme for \eqref{3.1} is presented.
The proposed scheme will be designed according to the algorithm in \cite{Ghazizadeh2019}:

\subsection{Finite-volume semi-discretization over quadtree grids}\label{S:3.1}
\begin{figure}[ht!]
\centerline{{\includegraphics[height=2.2in]{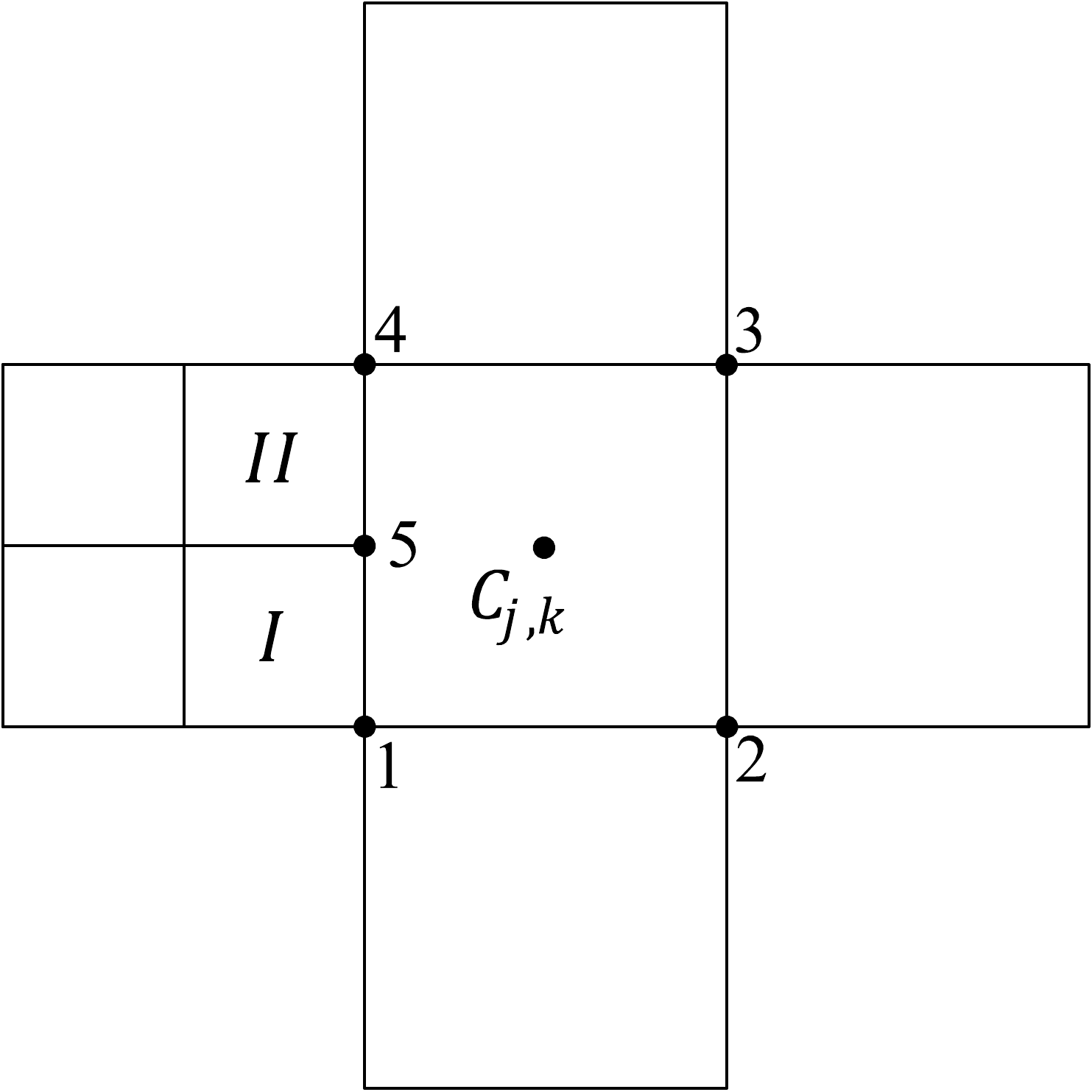}}}
\caption{\sf A configuration of cells neighboring $C_{j,k}$ in the regularised quadtree grid. \label{fig:2}}
\end{figure}

Let us represent each cell $C_{j,k}$ of size $\dx_{j,k}\times\dy_{j,k}$ centered at $(x_{j,k},y_{j,k})$ as 
a finite volume quadtree cell in the proposed scheme. The approximate averages of the cell read as:
\begin{equation}
\xbar{\bm{U}}_{j,k}(t)\approx
\frac{1}{\dx_{j,k}\dy_{j,k}}\int\limits_{x_\jmh}^{x_\jph}\int\limits_{y_\kmh}^{y_\kph}\bm{U}(x,y,t)\,{\rm d}y\,{\rm d}x,
\label{3.5}          
\end{equation}
where $x_{j\pm\hf}:=x_{j,k}\pm\dx_{j,k}/2$ and $y_{k\pm\hf}:=y_{j,k}\pm\dy_{j,k}/2$.
We present the proposed scheme for the configuration in Figure \ref{fig:2}. The left-neighboring cells of $C_{j,k}$ are denoted by $\Romannum{1}$ and
$\Romannum{2}$ which are centered at $(x_{j,k}-3\dx_{j,k}/4,y_{j,k}\pm\dy_{j,k}/4)$ with a size of $\dx_{j,k}/2\times\dy_{j,k}/2$.

The evolution of time-dependant cell averages $\,\xbar{\bm{U}}_{j,k}$, which are obtained after the semi-discretization of the system 
\eqref{3.1}--\eqref{3.4}, are computed by solving the following system of Ordinary Differential Equations (ODE):
\begin{equation}
\frac{{\rm d}}{{\rm d}t}\,\xbar{\bm{U}}_{j,k}=
-\frac{\bm{H}_{\jph,k}^x-\dfrac{\bm{H}_{\jmh,k-\frac{1}{4}}^x+\bm{H}_{\jmh,k+\frac{1}{4}}^x}{2}}{\dx_{j,k}}
-\frac{\bm{H}_{j,\kph}^y-\bm{H}_{j,\kmh}^y}{\dy_{j,k}}+\,\xbar{\bm{S}}_{j,k}.
\label{3.6}          
\end{equation}

In \eqref{3.6}, $\bm{H}_{\jph,k}^x$,
$\bm{H}_{\jmh,k\pm\frac{1}{4}}^x$, $\bm{H}_{j,\kph}^y$ and $\bm{H}_{j,\kmh}^y$ are the numerical fluxes, 
and $\xbar{\bm{S}}_{j,k}$ is a cell average of the source term:
\begin{equation}
\xbar{\bm{S}}_{j,k}\approx
\frac{1}{\dx_{j,k}\dy_{j,k}}\int\limits_{x_\jmh}^{x_\jph}\int\limits_{y_\kmh}^{y_\kph}\bm{S}(\bm{U},B)\,{\rm d}y\,{\rm d}x.
\label{3.7}          
\end{equation}

For the sake of brevity, we have omitted all of the time-dependent indexed quantities in \eqref{3.6}--\eqref{3.7}.

\subsection{Piecewise bilinear reconstruction of $B$}\label{S:3.5}
Cells of different sizes $\dx\times\dy$,
$\frac{\dx}{2}\times\frac{\dy}{2},\ldots,\frac{\dx}{2^{m-1}}\times\frac{\dy}{2^{m-1}}$ exist in the quadtree grid. 
The set of cells of the corresponding size is indicated 
by ${\cal C}^{(\ell)}$, that is, ${\cal C}^{(\ell)}=\{C_{j,k}: |C_{j,k}|=\frac{\dx}{2^{\ell-1}}\times\frac{\dy}{2^{\ell-1}}\}$. We exactly follow the steps that were introduced in \cite{Ghazizadeh2019} to reconstruct the bottom topography
$\widetilde B(x,y)$.

\subsection{Piecewise linear reconstructions}\label{S:3.4}
In this section, we construct a spatial second-order scheme, which employs a piecewise polynomial interpolation $\widetilde{\bm{U}}$, where
\begin{equation}
\widetilde{\bm{U}}(x,y)=(\bm{U}_x)_{j,k}[x-x_j]+(\bm{U}_y)_{j,k}[y-y_k],\quad(x,y)\in C_{j,k}.
\label{3.11}
\end{equation}

Such a reconstruction makes it impossible to develop a well-balanced scheme. 
Thus, instead of reconstructing conservative variables 
in $\bm{U}$, we reform $\bm{\cU}:=(w,hu,hv,\rho)^\top$ and then obtain the point values of 
$\bm{\cU}$ for the cell $C_{j,k}$ in Figure \ref{fig:2} which results in
\begin{equation}
\begin{aligned}
&\bm{\cU}_{\jph,k}^{+}=\,\xbar{\bm{\cU}}_{j+1,k}-\frac{\dx_{j+1,k}}{2}(\bm{\cU}_x)_{j+1,k}\,,\quad
\bm{\cU}_{\jph,k}^{-}=\,\xbar{\bm{\cU}}_{j,k}+\frac{\dx_{j,k}}{2}(\bm{\cU}_x)_{j,k}\,,\\
&\bm{\cU}_{\jmh,k\pm\frac{1}{4}}^{+}=\,\xbar{\bm{\cU}}_{j,k}-\frac{\dx_{j,k}}{2}(\bm{\cU}_x)_{j,k}\pm\frac{\dy_{j,k}}{2}(\bm{\cU}_y)_{j,k}\,,\\
&\bm{\cU}_{\jmh,k\pm\frac{1}{4}}^{-}=\,\xbar{\bm{\cU}}_{j-\frac{1}{4},k\pm\frac{1}{4}}+\frac{\dx_{j,k}}{4}(\bm{\cU}_x)_{j-\frac{1}{4},k\pm\frac{1}{4}},
\end{aligned}
\label{3.12} 
\end{equation}
where $\,\xbar{\bm{\cU}}$ denote the 
cell averages of $\bm{\cU}$. Note that in \eqref{3.12} the density variable $\rho_{j,k}$ is computed as
\begin{equation}
\rho_{j,k} := \frac{\xbar{(h\rho)}_{j,k}}{\xbar h_{j,k}}, \,\,\quad \xbar h_{j,k}:=\xbar w_{j,k}- B_{j,k},
\label{3.13} 
\end{equation}
and the point values of $h$ and $h\rho$ read as
$$
\begin{aligned}
&h_{\jmh,k\pm\frac{1}{4}}^+=w_{\jmh,k\pm\frac{1}{4}}^+-B_{\jmh,k\pm\frac{1}{4}},\quad h_{\jph,k}^-=w_{\jph,k}^--B_{\jph,k},\\
&h_{j,\kmh}^+=w_{j,\kmh}^+-B_{j,\kmh},\quad\mbox{and}\quad h_{j,\kph}^-=w_{j,\kph}^--B_{j,\kph},
\end{aligned}
$$
and
$$
\begin{aligned}
&(h\rho)_{\jmh,k\pm\frac{1}{4}}^+=h_{\jmh,k\pm\frac{1}{4}}^+\cdot\rho_{\jmh,k\pm\frac{1}{4}}^+,\quad (h\rho)_{\jph,k}^-=h_{\jph,k}^-\cdot\rho_{\jph,k}^-,\\
&(h\rho)_{j,\kmh}^+=h_{j,\kmh}^+\cdot\rho_{j,\kmh}^+,\quad\mbox{and}\quad (h\rho)_{j,\kph}^-=h_{j,\kph}^-\cdot\rho_{j,\kph}^-.
\end{aligned}
$$

We compute the slopes $(\bm{\cU}_x)$ and $(\bm{\cU}_y)$ by using the minmod limiter in order to minimize oscillations: 
\begin{equation}
\begin{aligned}
&(\bm{\cU}_x)_{j,k}={\rm minmod}\left(\frac{\,\xbar{\bm{\cU}}_{j,k}-\,\xbar{\bm{\cU}}_{j-\frac{1}{4},k-\frac{1}{4}}}{3\dx_{j,k}/4},\,
\frac{\,\xbar{\bm{\cU}}_{j,k}-\,\xbar{\bm{\cU}}_{j-\frac{1}{4},k+\frac{1}{4}}}{3\dx_{j,k}/4},\,
\frac{\,\xbar{\bm{\cU}}_{j+1,k}-\,\xbar{\bm{\cU}}_{j,k}}{\dx_{j,k}}\right),\\
&(\bm{\cU}_y)_{j,k}={\rm minmod}\left(\frac{\,\xbar{\bm{\cU}}_{j,k}-\,\xbar{\bm{\cU}}_{j,k-1}}{\dy_{j,k}},\,
\frac{\,\xbar{\bm{\cU}}_{j,k+1}-\,\xbar{\bm{\cU}}_{j,k}}{\dy_{j,k}}\right),
\end{aligned}
\label{3.14} 
\end{equation}
where the minmod function is defined by
\begin{equation*}
{\rm min}\{z_1,z_2,...\}:=\left\{
\begin{aligned}
&{\rm min}_j\{z_j\},&&\text{if}~z_j>0\quad\forall{j},\\
&{\rm max}_j\{z_j\},&&\text{if}~z_j<0\quad\forall{j},\\
&0,&&\text{otherwise}.
\end{aligned}\right.
\end{equation*}

Employing the minmod limiter \eqref{3.14} guarantees the positivity of the point values of $\rho$ \cite{Chertock2014}. For the
positivity-preserving correction of $w$, we exactly follow the steps in \cite{Ghazizadeh2019}. To prevent very small or even zero values of cell averages $\xbar h$, and point values of $h$ of cell $C_{j,k}$, the corresponding point values of $u$, $v$, and $\rho$ are computed
as follows \cite{Chertock2014, Ghazizadeh2019,Kurganov2007}:
\begin{equation*}
u:=\frac{\sqrt{2}\,h(hu)}{\sqrt{h^4+\max\{h^4,\varepsilon\}}},\qquad v:=\frac{\sqrt{2}\,h(hv)}{\sqrt{h^4+\max\{h^4,\varepsilon\}}}
,\qquad \rho:=\frac{\sqrt{2} \, h({h\rho})}{\sqrt{h^4+\max\{h^4,\varepsilon\}}},
\end{equation*}
where we choose $\varepsilon=\max\{\min_{j,k}\{(\dx_{j,k})^4\},\min_{j,k}\{(\dy_{j,k})^4\}\}$. 
The conservative variables recalculation is done by setting:
\begin{equation*}
(hu):=h\cdot u,\qquad (hv):=h\cdot v, \qquad (h\rho):=h\cdot \rho.
\end{equation*}

Note that all of the indices have been omitted in the above equations. 

\subsection{Local speeds}\label{S:3.3}
The one-sided local speeds of propagation, denoted at the corresponding cell interfaces by $a_{\alpha,\beta}^\pm$ and
$b_{\gamma,\delta}^\pm$ can be estimated by:
\begin{equation}
\begin{aligned}
&a_{\alpha,\beta}^+=\max\left\{u_{\alpha,\beta}^++\sqrt{\dfrac{g}{\rho_\circ} \,h_{\alpha,\beta}^+ \,\rho_{\alpha,\beta}^+},
\,u_{\alpha,\beta}^-+\sqrt{\dfrac{g}{\rho_\circ}\, h_{\alpha,\beta}^- \,\rho_{\alpha,\beta}^-},\,0\right\},\\
&a_{\alpha,\beta}^-=\min\left\{u_{\alpha,\beta}^+-\sqrt{\dfrac{g}{\rho_\circ} \, h_{\alpha,\beta}^+ \,\rho_{\alpha,\beta}^+},
\,u_{\alpha,\beta}^--\sqrt{\dfrac{g}{\rho_\circ} \, h_{\alpha,\beta}^- \,\rho_{\alpha,\beta}^-},\,0\right\},\\
&b_{\gamma,\delta}^+=\max\left\{v_{\gamma,\delta}^++\sqrt{\dfrac{g}{\rho_\circ} \, h_{\gamma,\delta}^+ \,\rho_{\gamma,\delta}^+},
\,v_{\gamma,\delta}^-+\sqrt{\dfrac{g}{\rho_\circ} \, h_{\gamma,\delta}^- \,\rho_{\gamma,\delta}^-},\,
0\right\},\\
&b_{\gamma,\delta}^-=\min\left\{v_{\gamma,\delta}^+-\sqrt{\dfrac{g}{\rho_\circ} \, h_{\gamma,\delta}^+ \,\rho_{\gamma,\delta}^+},
\,v_{\gamma,\delta}^--\sqrt{\dfrac{g}{\rho_\circ} \, h_{\gamma,\delta}^- \,\rho_{\gamma,\delta}^-},\,
0\right\}.
\end{aligned}
\label{3.17}
\end{equation}
where $(\alpha,\beta)\in\big\{(\jmh,k-\frac{1}{4}),(\jmh,k+\frac{1}{4}),(\jph,k)\big\}$ and
$(\gamma,\delta)\in\big\{(j,\kmh),(j,\kph)\big\}$ in Figure \ref{fig:2}.

\subsection{Central-upwind numerical fluxes}\label{S:3.2}
We use the central-upwind fluxes from \cite{Kurganov2007}:
\begin{equation}
\begin{aligned}
\bm{H}_{\alpha,\beta}^x=&\frac{a_{\alpha,\beta}^+\bm{F}(\bm{U}_{\alpha,\beta}^{-},B_{\alpha,\beta})-
a_{\alpha,\beta}^-\bm{F}(\bm{U}_{\alpha,\beta}^{+},B_{\alpha,\beta})}{a_{\alpha,\beta}^+-a_{\alpha,\beta}^-}
+\frac{a_{\alpha,\beta}^+a_{\alpha,\beta}^-}{a_{\alpha,\beta}^+-a_{\alpha,\beta}^-}\left[\bm{U}_{\alpha,\beta}^{+}-\bm{U}_{\alpha,\beta}^{-}
\right],\\
\bm{H}_{\gamma,\delta}^y=&\frac{b_{\gamma,\delta}^+\bm{G}(\bm{U}_{\gamma,\delta}^{-},B_{\gamma,\delta})-
b_{\gamma,\delta}^-\bm{G}(\bm{U}_{\gamma,\delta}^{+},B_{\gamma,\delta})}{b_{\gamma,\delta}^+-b_{\gamma,\delta}^-}
+\frac{b_{\gamma,\delta}^+b_{\gamma,\delta}^-}{b_{\gamma,\delta}^+-b_{\gamma,\delta}^-}
\left[\bm{U}_{\gamma,\delta}^{+}-\bm{U}_{\gamma,\delta}^{-}\right].
\end{aligned}
\label{3.18} 
\end{equation}

\subsection{Well-balanced discretization of the source term}\label{S:3.8}
When the discretized cell average of the source term,
$\,\xbar{\bm{S}}_{j,k}=\big(0,\,\xbar S_{j,k}^{\,(2)},\,\xbar S_{j,k}^{\,(3)}, 0\big)^\top$, exactly balances the numerical fluxes in Equation
\eqref{3.6} at the ``lake-at-rest'' steady state \eqref{1.2}, the numerical scheme is well-balanced. This means that the right-hand side (RHS) of \eqref{3.6} vanishes as long as
$\,\xbar{\bm{\cU}}_{j,k}\equiv\left(\widehat w,0,0, \widehat \rho\right)^\top$ for all $(j,k)$, where $\widehat w$ and $\widehat \rho$ are constants.

Notice that at the ``lake-at-rest'' state, all of the reconstructed point values are $w^\pm=\widetilde w$, $u^\pm=v^\pm=0$ and  $\rho^\pm=\widetilde \rho$, and thus,
$a^+_{\alpha,\beta}=-a^-_{\alpha,\beta},\,\forall(\alpha,\beta)$, $\,b^+_{\gamma,\delta}=-b^-_{\gamma,\delta}\,\forall(\gamma,\delta)$, and
the numerical fluxes \eqref{3.18} reduce to:
\begin{equation*}
\bm{H}_{\alpha,\beta}^x=\left(0,\dfrac{g}{2\rho_\circ}\widehat\rho\left(\widehat w-B_{\alpha,\beta}\right)^2, 0, 0\right)^\top,\quad
\bm{H}_{\gamma,\delta}^y=\left(0, 0,\dfrac{g}{2\rho_\circ}\widehat\rho\left(\widehat w-B_{\gamma,\delta}\right)^2, 0 \right)^\top,
\end{equation*}
and the flux terms on the RHS of \eqref{3.6} then become
\begin{equation}
\begin{aligned}
&-\frac{\bm{H}_{\jph,k}^x-\dfrac{\bm{H}_{\jmh,k-\frac{1}{4}}^x+\bm{H}_{\jmh,k+\frac{1}{4}}^x}{2}}{\dx_{j,k}}
-\frac{\bm{H}_{j,\kph}^y-\bm{H}_{j,\kmh}^y}{\dy_{j,k}}\\
&=-\dfrac{g}{2\rho_\circ}\widehat\rho\begin{pmatrix}
0\\
\dfrac{\left(\widehat w-B_{\jph,k}\right)^2}{\dx_{j,k}}-\dfrac{\left(\widehat w-B_{\jmh,k-\frac{1}{4}}\right)^2}{2\dx_{j,k}}-
\dfrac{\left(\widehat w-B_{\jmh,k+\frac{1}{4}}\right)^2}{2\dx_{j,k}}\\
\dfrac{\left(\widehat w-B_{j,\kph}\right)^2}{\dy_{j,k}}-\dfrac{\left(\widehat w-B_{j,\kmh}\right)^2}{\dy_{j,k}}
\end{pmatrix}.
\end{aligned}
\label{3.19}
\end{equation}

By applying Green's theorem, the source term in \ref{3.6} can be approximated by
$$
\begin{aligned}
&-\dfrac{g}{2\rho_\circ}\rho\,(w-B)B_x=\frac{g}{2\rho_\circ}\left[\rho\, (w-B)^2\right]_x-\frac{g}{\rho_\circ}\rho\, (w-B)w_x,\\
&-\dfrac{g}{2\rho_\circ}\rho\,(w-B)B_y=\frac{g}{2\rho_\circ}\left[\rho\, (w-B)^2\right]_y-\frac{g}{\rho_\circ}\rho\, (w-B)w_y.
\end{aligned}
$$

We now rewrite the cell averages of the second and third components of the integral in \eqref{3.7}
\begin{equation}
\dfrac{g}{2\rho_\circ}\int\limits_{y_\kmh}^{y_\kph}\left[\rho \,(w-B)^2\Big|_{x=x_\jph}-\rho \,(w-B)^2\Big|_{x=x_\jmh}\right]\,{\rm d}y-
\dfrac{g}{\rho_\circ}\int\limits_{x_\jmh}^{x_\jph}\int\limits_{y_\kmh}^{y_\kph}\rho\, (w-B)w_x\,{\rm d}y\,{\rm d}x,
\label{3.20}
\end{equation}
and
\begin{equation}
\dfrac{g}{2\rho_\circ}\int\limits_{x_\jmh}^{x_\jph}\left[\rho\,(w-B)^2\Big|_{y=y_\kph}-\rho\,(w-B)^2\Big|_{y=y_\kmh}\right]\,{\rm d}x-
\dfrac{g}{\rho_\circ}\int\limits_{x_\jmh}^{x_\jph}\int\limits_{y_\kmh}^{y_\kph}\rho\, (w-B)w_y\,{\rm d}y\,{\rm d}x,
\label{3.21}
\end{equation}

We then approximate the integrals in \eqref{3.20} and \eqref{3.21}, 
which results in the following quadrature for the second and third components of the
source term \cite{Ghazizadeh2019}:
\begin{equation}
\begin{aligned}
\xbar{S}_{j,k}^{\,(2)}\approx\frac{g}{2\rho_\circ \dx_{j,k}}\Bigg[& \rho_{\jph,k}^-\left(w_{\jph,k}^--B_{\jph,k}\right)^2-
\frac{\rho_{\jmh,k-\frac{1}{4}}^+\left(w_{\jmh,k-\frac{1}{4}}^+-B_{\jmh,k-\frac{1}{4}}\right)^2}{2}\\
&-\frac{\rho_{\jmh,k+\frac{1}{4}}^+\left(w_{\jmh,k+\frac{1}{4}}^+-B_{\jmh,k+\frac{1}{4}}\right)^2}{2}\Bigg]-
\dfrac{g}{\rho_\circ}\rho_{j,k}(w_x)_{j,k}\left(\,\xbar w_{j,k}-B_{j,k}\right),\\
\xbar{S}_{j,k}^{\,(3)}\approx\frac{g}{2\rho_\circ\dy_{j,k}}\bigg[&\rho_{j,\kph}^-\left(w_{j,\kph}^--B_{j,\kph}\right)^2-
\rho_{j,\kmh}^+\left(w_{j,\kmh}^+-B_{j,\kmh}\right)^2
\bigg]\\
&-\dfrac{g}{\rho_\circ}\rho_{j,k}(w_y)_{j,k}\left(\,\xbar w_{j,k}-B_{j,k}\right).
\end{aligned}
\label{3.22}
\end{equation}

We finally state that the scheme now preserves the solution at ``lake-at-rest'' where, $(w_x)_{j,k}=(w_y)_{j,k}\equiv0,\,\forall(j,k)$ 
and thus \eqref{3.19} and \eqref{3.22} 
express that the RHS of \eqref{3.6} vanishes and therefore, the scheme is well-balanced.

\subsection{Positivity-preserving property}\label{S:3.7}
In this section, we extend the positivity-preserving proof from \cite{Ghazizadeh2019} to implement on the coupled variable density system. 
We use a forward Euler method to integrate
Equation \eqref{3.6} in time, which results in
\begin{equation}
\xbar w_{j,k}^{\,n+1}=\,\xbar w_{j,k}^{\,n}-\lambda_{j,k}^n\left(H_{\jph,k}^{x,(1)}-\frac{H_{\jmh,k-\frac{1}{4}}^{x,(1)}+
H_{\jmh,k+\frac{1}{4}}^{x,(1)}}{2}\right)-\mu_{j,k}^n\left(H_{j,\kph}^{y,(1)}-H_{j,\kmh}^{y,(1)}\right),
\label{3.23}
\end{equation}
\begin{equation}
\xbar {(h\rho)}_{j,k}^{\,n+1}=\,\xbar {(h\rho)}_{j,k}^{\,n}-\lambda_{j,k}^n\left(H_{\jph,k}^{x,(4)}-\frac{H_{\jmh,k-\frac{1}{4}}^{x,(4)}+
H_{\jmh,k+\frac{1}{4}}^{x,(4)}}{2}\right)-\mu_{j,k}^n\left(H_{j,\kph}^{y,(4)}-H_{j,\kmh}^{y,(4)}\right),
\label{3.24}
\end{equation}
where $\,\xbar w_{j,k}^{\,n}:=\,\xbar w_{j,k}(t^n)$, $\,\xbar w_{j,k}^{\,n+1}:=\,\xbar w_{j,k}(t^{n+1})$, $\,
\xbar {(h\rho)}_{j,k}^{\,n}:=\,\xbar  {(h\rho)}_{j,k}(t^n)$, and $\,\xbar  {(h\rho)}_{j,k}^{\,n+1}:=\,\xbar  {(h\rho)}_{j,k}(t^{n+1})$ with $t^{n+1}=t^n+\dt^n$,
$\lambda_{j,k}^n:=\dt^n/\dx_{j,k}$, $\mu_{j,k}^n:=\dt^n/\dy_{j,k}$, and the numerical fluxes on the RHS are evaluated at time level $t=t^n$ 
using \eqref{3.18}:
\begin{equation}
\begin{aligned}
H_{\alpha,\beta}^{x,(1)}=&\frac{a_{\alpha,\beta}^+(hu)_{\alpha,\beta}^--a_{\alpha,\beta}^-(hu)_{\alpha,\beta}^+}
{a_{\alpha,\beta}^+-a_{\alpha,\beta}^-}+\frac{a_{\alpha,\beta}^+a_{\alpha,\beta}^-}{a_{\alpha,\beta}^+-a_{\alpha,\beta}^-}
\left[w_{\alpha,\beta}^+-w_{\alpha,\beta}^-\right],\\
H_{\gamma,\delta}^{y,(1)}=&\frac{b_{\gamma,\delta}^+(hv)_{\gamma,\delta}^--b_{\gamma,\delta}^-(hv)_{\gamma,\delta}^+}
{b_{\gamma,\delta}^+-b_{\gamma,\delta}^-}+\frac{b_{\gamma,\delta}^+b_{\gamma,\delta}^-}{b_{\gamma,\delta}^+-b_{\gamma,\delta}^-}
\left[w_{\gamma,\delta}^+-w_{\gamma,\delta}^-\right],
\end{aligned}
\label{3.25}
\end{equation}
and
\begin{equation}
\begin{aligned}
H_{\alpha,\beta}^{x,(4)}=&\frac{a_{\alpha,\beta}^+\,\rho_{\alpha,\beta}^-(hu)_{\alpha,\beta}^--a_{\alpha,\beta}^-\,\rho_{\alpha,\beta}^+(hu)_{\alpha,\beta}^+}
{a_{\alpha,\beta}^+-a_{\alpha,\beta}^-}+\frac{a_{\alpha,\beta}^+a_{\alpha,\beta}^-}{a_{\alpha,\beta}^+-a_{\alpha,\beta}^-}
\left[{(h\rho)}_{\alpha,\beta}^+-{(h\rho)}_{\alpha,\beta}^-\right],\\
H_{\gamma,\delta}^{y,(4)}=&\frac{b_{\gamma,\delta}^+\,\rho_{\gamma,\delta}^-(hv)_{\gamma,\delta}^--b_{\gamma,\delta}^-\,\rho_{\gamma,\delta}^+(hv)_{\gamma,\delta}^+}
{b_{\gamma,\delta}^+-b_{\gamma,\delta}^-}+\frac{b_{\gamma,\delta}^+b_{\gamma,\delta}^-}{b_{\gamma,\delta}^+-b_{\gamma,\delta}^-}
\left[{(h\rho)}_{\gamma,\delta}^+-{(h\rho)}_{\gamma,\delta}^-\right].
\end{aligned}
\label{3.26}
\end{equation}

If $\,\xbar h_{j,k}^{\,n}\ge0$ for all $(j,k)$, then the point values of the computed $h$ are nonnegative \cite{Kurganov2007}. Moreover, using the bilinear piece for the bottom topography from \cite{Ghazizadeh2019} and the similar
relationships for the reconstructed point values of $w$, we have
\begin{equation}
\xbar h_{j,k}^{\,n}=\frac{1}{4}\left(\frac{h_{\jmh,k-\frac{1}{4}}^++h_{\jmh,k+\frac{1}{4}}^+}{2}+h_{\jph,k}^-+h_{j,\kmh}^++h_{j,\kph}^-\right)
\label{3.27}
\end{equation}
for the grid configuration in Figure \ref{fig:2}.

We now subtract $B_{j,k}$ from both sides of \eqref{3.23} and use \eqref{3.25} and \eqref{3.27} to rewrite \eqref{3.23} as follows:
\begin{equation}
\begin{aligned}
\xbar h_{j,k}^{\,n+1}=&-\lambda_{j,k}^na_{\jph,k}^-\cdot\frac{a_{\jph,k}^+-u_{\jph,k}^+}{a_{\jph,k}^+-a_{\jph,k}^-}\cdot h_{\jph,k}^++
\left[\frac{1}{4}-\lambda_{j,k}^na_{\jph,k}^+\cdot\frac{u_{\jph,k}^--a_{\jph,k}^-}{a_{\jph,k}^+-a_{\jph,k}^-}\right]h_{\jph,k}^-\\
&+\frac{\lambda_{j,k}^na_{\jmh,k-\frac{1}{4}}^+}{2}\cdot\frac{u_{\jmh,k-\frac{1}{4}}^--a_{\jmh,k-\frac{1}{4}}^-}
{a_{\jmh,k-\frac{1}{4}}^+-a_{\jmh,k-\frac{1}{4}}^-}\cdot h_{\jmh,k-\frac{1}{4}}^+\\
&+\hf\left[\frac{1}{4}-\lambda_{j,k}^na_{\jmh,k-\frac{1}{4}}^-\cdot\frac{a_{\jmh,k-\frac{1}{4}}^+-u_{\jmh,k-\frac{1}{4}}^+}
{a_{\jmh,k-\frac{1}{4}}^+-a_{\jmh,k-\frac{1}{4}}^-}\right]h_{\jmh,k-\frac{1}{4}}^+\\
&+\frac{\lambda_{j,k}^na_{\jmh,k+\frac{1}{4}}^+}{2}\cdot\frac{u_{\jmh,k+\frac{1}{4}}^--a_{\jmh,k+\frac{1}{4}}^-}
{a_{\jmh,k+\frac{1}{4}}^+-a_{\jmh,k+\frac{1}{4}}^-}\cdot h_{\jmh,k+\frac{1}{4}}^+\\
&+\hf\left[\frac{1}{4}-\lambda_{j,k}^na_{\jmh,k+\frac{1}{4}}^-\cdot\frac{a_{\jmh,k+\frac{1}{4}}^+-u_{\jmh,k+\frac{1}{4}}^+}
{a_{\jmh,k+\frac{1}{4}}^+-a_{\jmh,k+\frac{1}{4}}^-}\right]h_{\jmh,k+\frac{1}{4}}^+\\
&-\mu_{j,k}^nb_{j,\kph}^-\cdot\frac{b_{j,\kph}^+-v_{j,\kph}^+}{b_{j,\kph}^+-b_{j,\kph}^-}\cdot h_{j,\kph}^++
\left[\frac{1}{4}-\mu_{j,k}^nb_{j,\kph}^+\cdot\frac{v_{j,\kph}^--b_{j,\kph}^-}{b_{j,\kph}^+-b_{j,\kph}^-}\right]h_{j,\kph}^-\\
&+\mu_{j,k}^nb_{j,\kmh}^+\cdot\frac{v_{j,\kmh}^--b_{j,\kmh}^-}{b_{j,\kmh}^+-b_{j,\kmh}^-}\cdot h_{j,\kmh}^-+
\left[\frac{1}{4}+\mu_{j,k}^nb_{j,\kmh}^-\cdot\frac{b_{j,\kmh}^+-v_{j,\kmh}^+}{b_{j,\kmh}^+-b_{j,\kmh}^-}\right]h_{j,\kmh}^+.
\end{aligned}
\label{3.28} 
\end{equation}

This shows that the cell averages of $h$ are linear combinations of the reconstructed nonnegative point
values of $h$. Thus, $\,\xbar h_{j,k}^{\,n+1}\ge0$ provided all of the coefficients in this linear combination are nonnegative.

One can obtain a similar proof for positivity of $\xbar {(h\rho)}_{j,k}^{\,n+1}$ by using the following statements \cite{Chertock2014},
\begin{equation}
\xbar{{h}}_{j,k}=\frac{{h}_{\jph,k}^{-}  + \cfrac{h_{\jmh,k-\frac{1}{4}}^{+} +  h_{\jmh,k+\frac{1}{4}}^{+}}{2}}{2}\,,\quad
{\rho}_{j,k}=\frac{{\rho}_{\jph,k}^{-}  + \cfrac{\rho_{\jmh,k-\frac{1}{4}}^{+} + \rho_{\jmh,k+\frac{1}{4}}^{+}}{2}}{2},
\label{3.29} 
\end{equation}
and thereby, utilizing \eqref{3.13}, one obtains:

\begin{equation}
\begin{aligned}
\xbar{{(h\rho)}}^n_{j,k}= \,&\frac{1}{4} \left[ h_{\jph,k}^{-} \,\rho_{\jph,k}^{-}  + \frac{h_{\jmh,k-\frac{1}{4}}^{+} \,
\rho_{\jmh,k-\frac{1}{4}}^{+} +  h_{\jmh,k+\frac{1}{4}}^{+} \,\rho_{\jmh,k+\frac{1}{4}}^{+}}{2}\right] \\
& + \frac{1}{4} \left[\frac{ h_{\jph,k}^{-} \,\rho_{\jmh,k-\frac{1}{4}}^{+} + h_{\jph,k}^{-} \,\rho_{\jmh,k+\frac{1}{4}}^{+}   +    
h_{\jmh,k-\frac{1}{4}}^{+} \,\rho_{\jph,k}^{-} +h_{\jmh,k+\frac{1}{4}}^{+} \,\rho_{\jph,k}^{-}}{2}   \right] \\
& + \frac{1}{4} \left[\frac{h_{\jmh,k-\frac{1}{4}}^{+} \,\rho_{\jmh,k+\frac{1}{4}}^{+} +  h_{\jmh,k+\frac{1}{4}}^{+} \,\rho_{\jmh,k-\frac{1}{4}}^{+}}{4}   \right].
\end{aligned}
\label{3.30} 
\end{equation}

Similarly, it can be shown that

\begin{equation}
\xbar{{(h\rho)}}^n_{j,k}= \,\frac{1}{4} \left[ h_{j,\kph}^{-} \,\rho_{j,\kph}^{-} +  h_{j,\kmh}^{+} \,\rho_{j,\kmh}^{+} \right] + 
\frac{1}{4} \left[ h_{j,\kph}^{-} \,\rho_{j,\kmh}^{+} +  h_{j,\kmh}^{+} \,\rho_{j,\kph}^{-} \right].
\label{3.31} 
\end{equation}

Finally, from \eqref{3.30} and \eqref{3.31} we have

\begin{equation}
\begin{aligned}
\xbar{{(h\rho)}}^n_{j,k}= \,&\frac{1}{8} \left[ h_{\jph,k}^{-} \,\rho_{\jph,k}^{-}  +  
h_{j,\kph}^{-} \,\rho_{j,\kph}^{-} +  h_{j,\kmh}^{+} \,\rho_{j,\kmh}^{+} \right] \\
& + \frac{1}{8} \left[\frac{h_{\jmh,k-\frac{1}{4}}^{+} \,\rho_{\jmh,k-\frac{1}{4}}^{+} +  
h_{\jmh,k+\frac{1}{4}}^{+} \,\rho_{\jmh,k+\frac{1}{4}}^{+}}{2}\right] \\
& + \frac{1}{8} \left[\frac{ h_{\jph,k}^{-} \,\rho_{\jmh,k-\frac{1}{4}}^{+} + h_{\jph,k}^{-} \,\rho_{\jmh,k+\frac{1}{4}}^{+}   +    
h_{\jmh,k-\frac{1}{4}}^{+} \,\rho_{\jph,k}^{-} +h_{\jmh,k+\frac{1}{4}}^{+} \,\rho_{\jph,k}^{-}}{2}   \right] \\
& + \frac{1}{8} \left[\frac{h_{\jmh,k-\frac{1}{4}}^{+} \,\rho_{\jmh,k+\frac{1}{4}}^{+} +  
h_{\jmh,k+\frac{1}{4}}^{+} \,\rho_{\jmh,k-\frac{1}{4}}^{+}}{4}   \right].
\end{aligned}
\label{3.32} 
\end{equation}

We rewrite \eqref{3.24} as follows:

\begin{equation}
\begin{aligned}
\xbar {(h\rho)}_{j,k}^{\,n+1}=&-\lambda_{j,k}^na_{\jph,k}^-\cdot\frac{a_{\jph,k}^+-u_{\jph,k}^+}{a_{\jph,k}^+-a_{\jph,k}^-}\cdot h_{\jph,k}^+ \,\rho_{\jph,k}^+ \\
&+\left[\frac{1}{8}-\lambda_{j,k}^na_{\jph,k}^+\cdot\frac{u_{\jph,k}^--a_{\jph,k}^-}{a_{\jph,k}^+-a_{\jph,k}^-}\right]h_{\jph,k}^-\,\rho_{\jph,k}^-\\
&+\frac{\lambda_{j,k}^na_{\jmh,k-\frac{1}{4}}^+}{2}\cdot\frac{u_{\jmh,k-\frac{1}{4}}^--a_{\jmh,k-\frac{1}{4}}^-}
{a_{\jmh,k-\frac{1}{4}}^+-a_{\jmh,k-\frac{1}{4}}^-}\cdot h_{\jmh,k-\frac{1}{4}}^+  \,\rho_{\jmh,k-\frac{1}{4}}^+\\
&+\hf\left[\frac{1}{8}-\lambda_{j,k}^na_{\jmh,k-\frac{1}{4}}^-\cdot\frac{a_{\jmh,k-\frac{1}{4}}^+-u_{\jmh,k-\frac{1}{4}}^+}
{a_{\jmh,k-\frac{1}{4}}^+-a_{\jmh,k-\frac{1}{4}}^-}\right]h_{\jmh,k-\frac{1}{4}}^+ \,\rho_{\jmh,k-\frac{1}{4}}^+\\
&+\frac{\lambda_{j,k}^na_{\jmh,k+\frac{1}{4}}^+}{2}\cdot\frac{u_{\jmh,k+\frac{1}{4}}^--a_{\jmh,k+\frac{1}{4}}^-}
{a_{\jmh,k+\frac{1}{4}}^+-a_{\jmh,k+\frac{1}{4}}^-}\cdot h_{\jmh,k+\frac{1}{4}}^+  \,\rho_{\jmh,k+\frac{1}{4}}^+\\
&+\hf\left[\frac{1}{8}-\lambda_{j,k}^na_{\jmh,k+\frac{1}{4}}^-\cdot\frac{a_{\jmh,k+\frac{1}{4}}^+-u_{\jmh,k+\frac{1}{4}}^+}
{a_{\jmh,k+\frac{1}{4}}^+-a_{\jmh,k+\frac{1}{4}}^-}\right]h_{\jmh,k+\frac{1}{4}}^+ \,\rho_{\jmh,k+\frac{1}{4}}^+\\
&-\mu_{j,k}^nb_{j,\kph}^-\cdot\frac{b_{j,\kph}^+-v_{j,\kph}^+}{b_{j,\kph}^+-b_{j,\kph}^-}\cdot h_{j,\kph}^+ \,\rho_{j,\kph}^+\\
&+\left[\frac{1}{8}-\mu_{j,k}^nb_{j,\kph}^+\cdot\frac{v_{j,\kph}^--b_{j,\kph}^-}{b_{j,\kph}^+-b_{j,\kph}^-}\right]h_{j,\kph}^- \,\rho_{j,\kph}^-\\
&+\mu_{j,k}^nb_{j,\kmh}^+\cdot\frac{v_{j,\kmh}^--b_{j,\kmh}^-}{b_{j,\kmh}^+-b_{j,\kmh}^-}\cdot h_{j,\kmh}^- \, \rho_{j,\kmh}^-\\
&+\left[\frac{1}{8}+\mu_{j,k}^nb_{j,\kmh}^-\cdot\frac{b_{j,\kmh}^+-v_{j,\kmh}^+}{b_{j,\kmh}^+-b_{j,\kmh}^-}\right]h_{j,\kmh}^+ \, \rho_{j,\kmh}^+\\
& + \frac{1}{8} \left[\frac{ h_{\jph,k}^{-} \,\rho_{\jmh,k-\frac{1}{4}}^{+} + h_{\jph,k}^{-} \,\rho_{\jmh,k+\frac{1}{4}}^{+}   +   
 h_{\jmh,k-\frac{1}{4}}^{+} \,\rho_{\jph,k}^{-} +h_{\jmh,k+\frac{1}{4}}^{+} \,\rho_{\jph,k}^{-}}{2}   \right] \\
& + \frac{1}{8} \left[\frac{h_{\jmh,k-\frac{1}{4}}^{+} \,\rho_{\jmh,k+\frac{1}{4}}^{+} +  h_{\jmh,k+\frac{1}{4}}^{+} \,\rho_{\jmh,k-\frac{1}{4}}^{+}}{4}   \right].
\end{aligned}
\label{3.33} 
\end{equation}

From \eqref{3.12} and positivity-preserving correction of $w$, which guarantee the positivity of the reconstructed point values of $h$, the last two terms in \eqref{3.33} are nonnegative. The other terms are also nonnegative similar to those in \eqref{3.28}. Note that using the minmod limiter \eqref{3.14} guarantees positivity of the point values of $\rho$.
Since both $h^{n+1}_{j,k}$ and $\xbar {(h\rho)}_{j,k}^{\,n+1}$ are nonnegative, thus $\rho^{n+1}_{j,k} = \xbar {(h\rho)}_{j,k}^{\,n+1} / h^{n+1}_{j,k}$ 
is nonnegative as well.


We use \eqref{3.17} to satisfy the Courant-Friedrichs-Lewy (CFL) type condition as follows:
$$
\dt\le\frac{1}{8}\min\left[\min_{j,k}\left\{
\frac{\dx_{j,k}}{\max\limits_{(\alpha,\beta)}\left[\max\left\{a_{\alpha,\beta}^+,-a_{\alpha,\beta}^-\right\}\right]}\right\},\,
\min_{j,k}\left\{
\frac{\dy_{j,k}}{\max\limits_{(\gamma,\delta)}\left[\max\left\{b_{\gamma,\delta}^+,-b_{\gamma,\delta}^-\right\}\right]}\right\}\right].
$$

Finally, in all of the numerical experiments, we have used the
three-stage third-order strong stability preserving (SSP) Runge-Kutta solver (see, e.g., \cite{Ghazizadeh2019,GKS,GST}). 

\subsection{Quadtree grid adaptivity}\label{S:3.10}
At the new time level $t=t^{n+1}$, the quadtree grid locally refines or coarsens for the 
next timestep. We first need to compute the slopes $\{(w_x)_{j,k}^{n+1}\}$ and $\{(w_y)_{j,k}^{n+1}\}$, and  
$\{(\rho_x)_{j,k}^{n+1}\}$ and $\{(\rho_y)_{j,k}^{n+1}\}$  on the old grid
(which is denoted by $\{C_{j,k}^{\rm old}\}$) according to \S\ref{S:3.4} and then select the centers of those cells $C_{j,k}^{\rm old}$, at
which (see \cite{Ghazizadeh2019}):
\begin{equation}
(w_x)_{j,k}^{n+1}\ge C_{w,\,\rm seed}\quad\mbox{or}\quad(w_y)_{j,k}^{n+1}\ge C_{w,\,\rm seed},
\label{3.34}
\end{equation}
and
\begin{equation}
(\rho_x)_{j,k}^{n+1}\ge C_{\rho,\,\rm seed}\quad\mbox{or}\quad(\rho_y)_{j,k}^{n+1}\ge C_{\rho,\,\rm seed}.
\label{3.35}
\end{equation}

We denote the required seeding points to generate the new grid by $\{C_{j,k}^{\rm new}\}$. In \eqref{3.34} and \eqref{3.35}, 
$C_{w, \,\rm seed}$ and $C_{\rho, \,\rm seed}$ are
constants that depend on the problem at hand (e.g. the maximum level of the quadtree,
the Froude number, and the bottom topography function). When the grid locally refines or coarsens, at the end of the evolution step, the solution in terms of the computed cell
averages $\big\{\big(\,\xbar{\bm{\cU}}_{j,k}^{\,n+1}\big)_{\rm old}\big\}$ over the grid $\{C_{j,k}^{\rm old}\}$, should be projected onto
the new grid $\{C_{j,k}^{\rm new}\}$. This should be done in a conservative manner as follows:

{\em Case 1:} If $C_{j,k}^{\rm new}=C_{j',k'}^{\rm old}$ for some $(j',k')$, that is, if the cell $C_{j',k'}^{\rm old}$ does not need to be
refined/coarsened, then
\begin{equation*}
\big(\,\xbar{\bm{\cU}}_{j,k}^{\,n+1}\big)_{\rm new}=\big(\,\xbar{\bm{\cU}}_{j',k'}^{\,n+1}\big)_{\rm old}.
\end{equation*}

{\em Case 2:} If $C_{j,k}^{\rm new}\in{\cal C}^{\ell+p}$ is a ``child'' cell of $C_{j',k'}^{\rm old}\in{\cal C}^{\ell}$ for some $j'$, $k'$
and $p>0$ (that is, if the cell $C_{j',k'}^{\rm old}$ was refined and $C_{j,k}^{\rm new}\subset C_{j',k'}^{\rm old}$), then
\begin{equation*}
\big(\,\xbar{\bm{\cU}}^{\,n+1}_{j,k}\big)_{\rm new}=\big(\,\xbar{\bm{\cU}}^{\,n+1}_{j',k'}\big)_{\rm old}+
\big((\bm{\cU}_x)^{n+1}_{j',k'}\big)_{\rm old}\left[x_j^{\rm new}-x_{j'}^{\rm old}\right]+
\big((\bm{\cU}_y)^{n+1}_{j',k'}\big)_{\rm old}\left[y_k^{\rm new}-y_{k'}^{\rm old}\right].
\end{equation*}

{\em Case 3:} If $C_{j,k}^{\rm new}\in{\cal C}^{\ell-p}$ is a ``parent'' cell of $C_{j',k'}^{\rm old}\in{\cal C}^{\ell}$ for some $j'$, $k'$
and $p>0$ (that is, if the cell $C_{j',k'}^{\rm old}$ was coarsened and $C_{j,k}^{\rm new}\supset C_{j',k'}^{\rm old}$), then
\begin{equation*}
\big(\,\xbar{\bm{\cU}}^{\,n+1}_{j,k}\big)_{\rm new}=\frac{1}{4^p}\underset{j'',k'':\,C_{j'',k''}^{\rm old}\subset C_{j,k}^{\rm new}}{\sum\sum}
\big(\,\xbar{\bm{\cU}}^{\,n+1}_{j'',k''}\big)_{\rm old}.
\end{equation*}

\section{Numerical examples}\label{S:4}
In this section, we investigate the performance of the proposed scheme in four numerical examples. In all of the examples, we take
$g=1$ and $\rho_\circ = 997$.

\subsection*{Example 1 --- Circular dam break with constant density}
We demonstrate the ability of the proposed scheme to generate adaptive grids at each timestep and 
maintain symmetry in this example. A circular water column collapses on a horizontal flat bottom topography $40 \times 40$ dimensions plane
\cite{EleuterioFToro2001, Jiang2011} where:
\begin{align*}
&w(x,y,0)=\left\{\begin{array}{lc}2,&~(x-20)^2+(y-20)^2< 2.5^2,\\1,&\mbox{otherwise},\end{array}\right.\\
&u(x,y,0)=v(x,y,0)\equiv0, \quad \rho(x,y,0) \equiv \rho_\circ.
\end{align*}

\begin{figure}[ht!]
\centerline{\includegraphics[height=3.8cm]{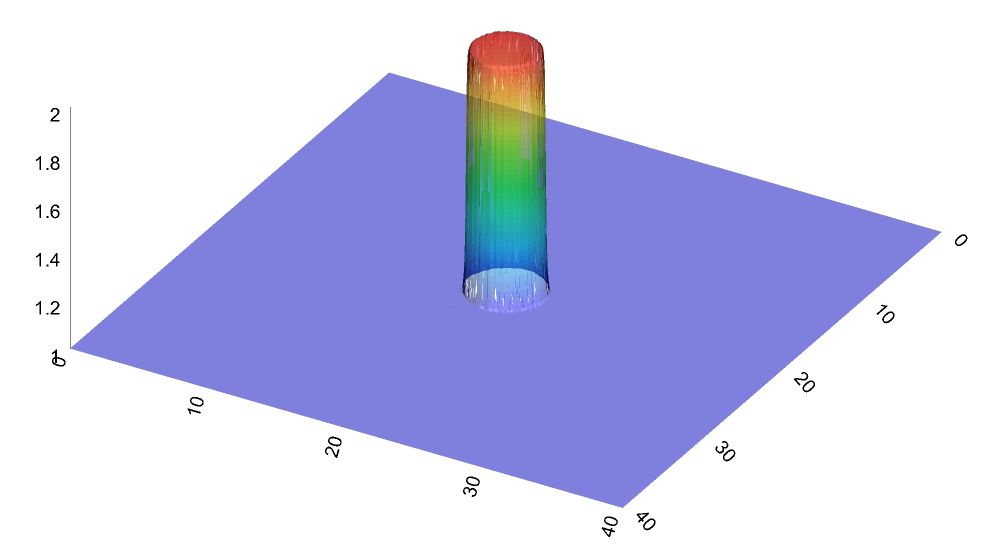}\hspace*{1cm}\includegraphics[height=3.8cm]{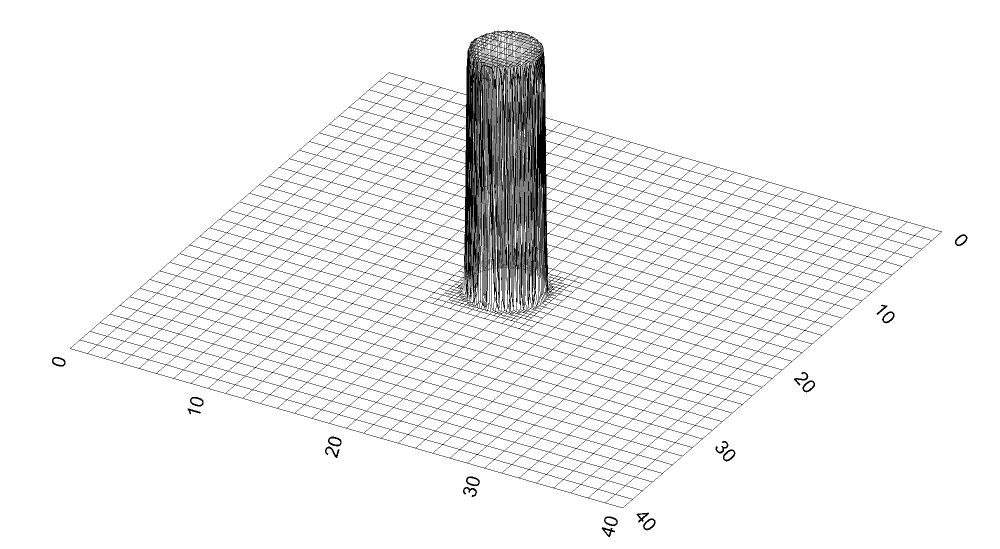}}
\vspace*{0.25cm}
\centerline{\includegraphics[height=3.8cm]{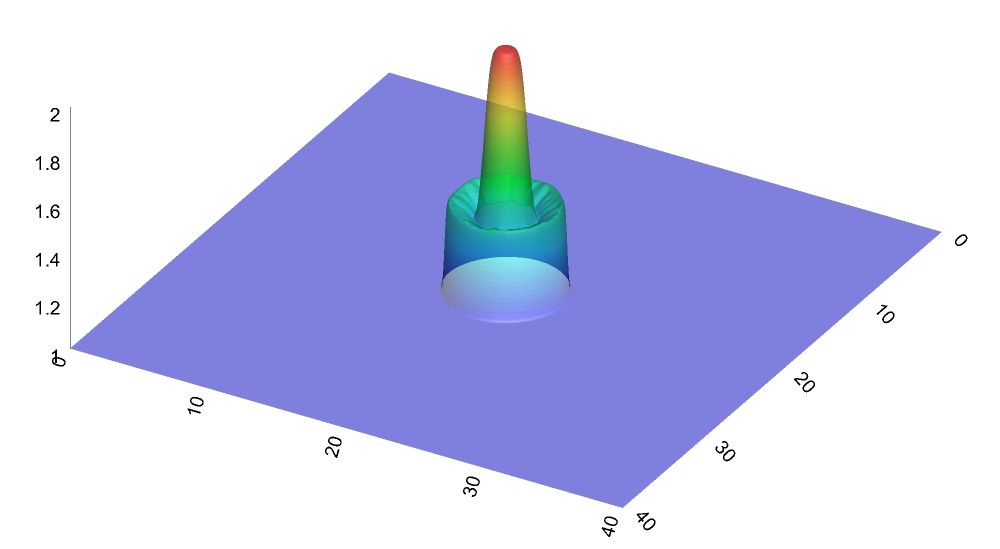}\hspace*{1cm}\includegraphics[height=3.8cm]{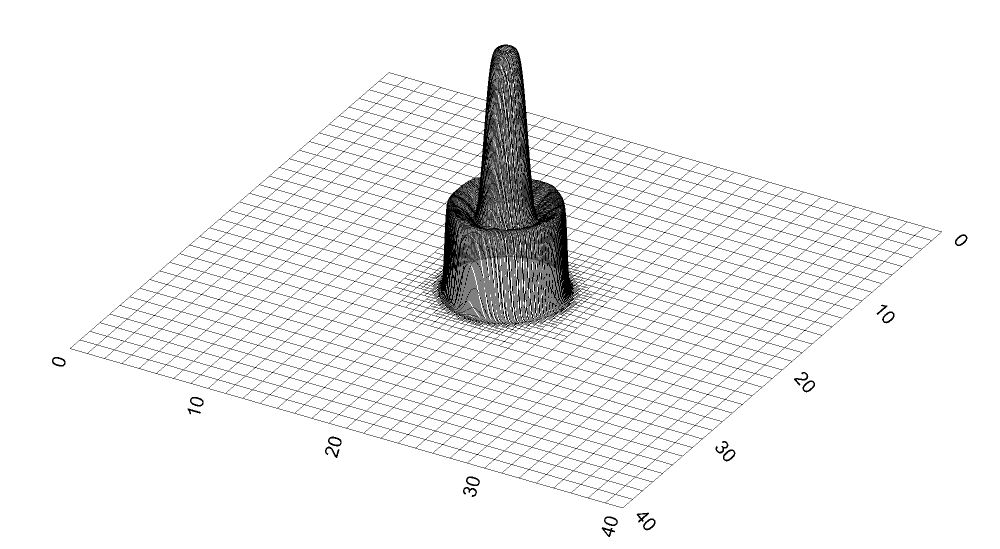}}
\vspace*{0.25cm}
\centerline{\includegraphics[height=3.8cm]{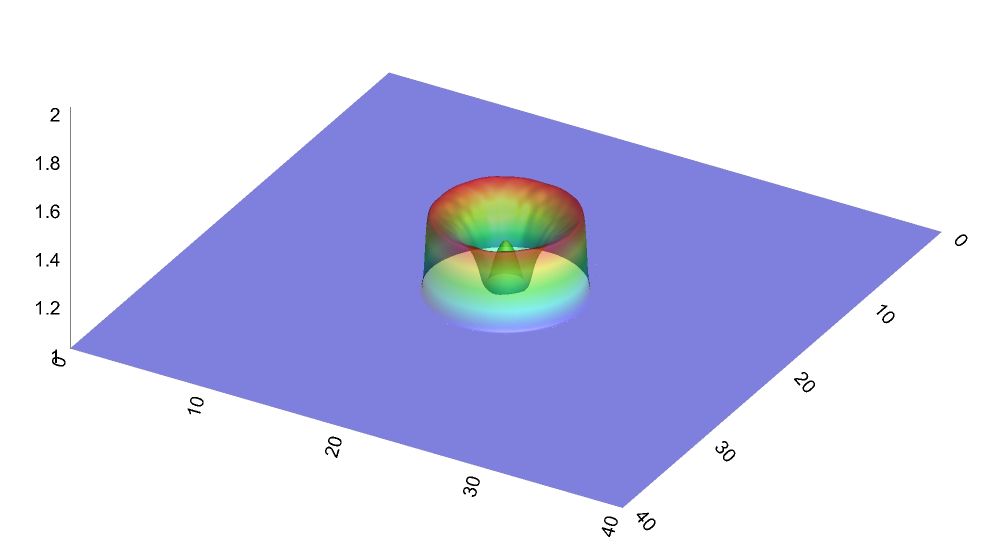}\hspace*{1cm}\includegraphics[height=3.8cm]{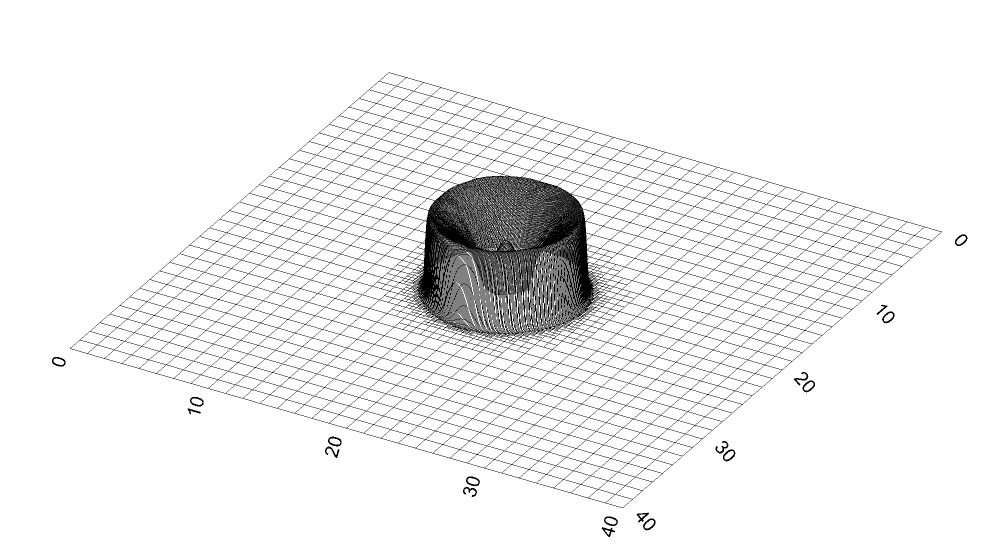}}
\vspace*{0.25cm}
\centerline{\includegraphics[height=3.8cm]{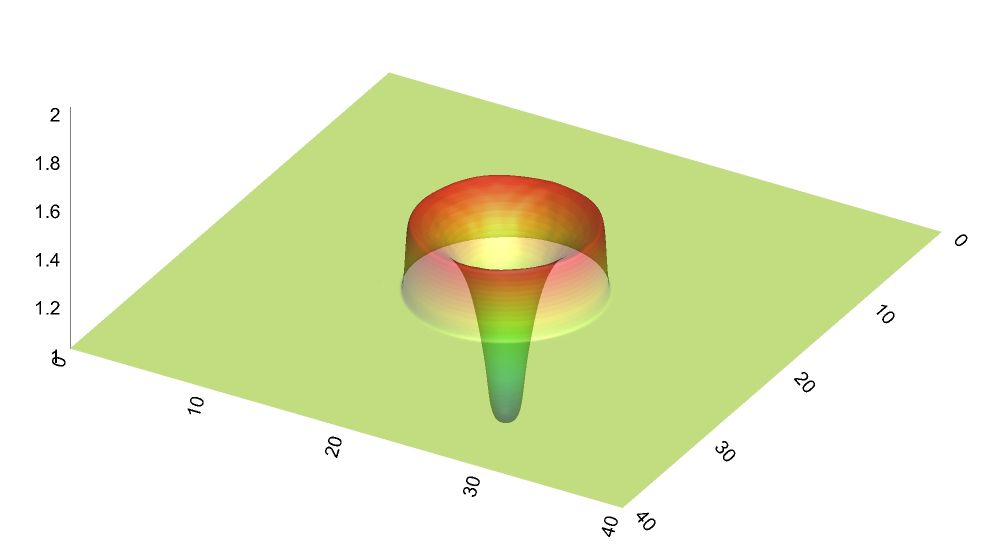}\hspace*{1cm}\includegraphics[height=3.8cm]{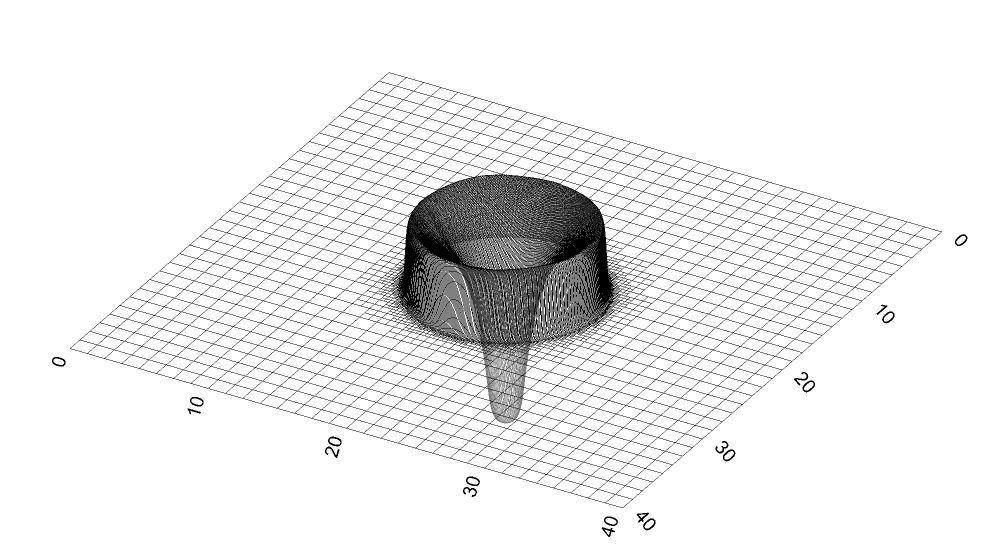}}
\vspace*{0.25cm}
\centerline{\includegraphics[height=3.8cm]{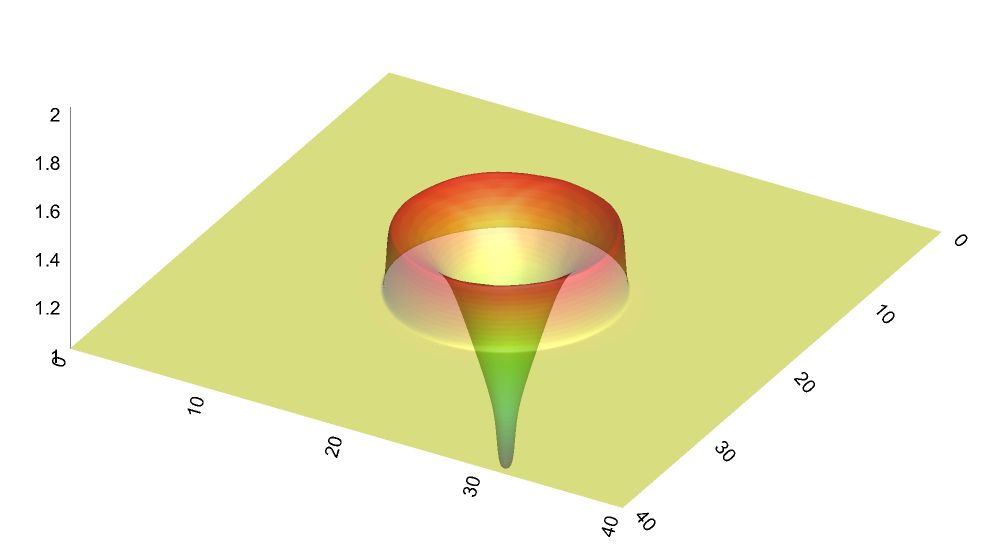}\hspace*{1cm}\includegraphics[height=3.8cm]{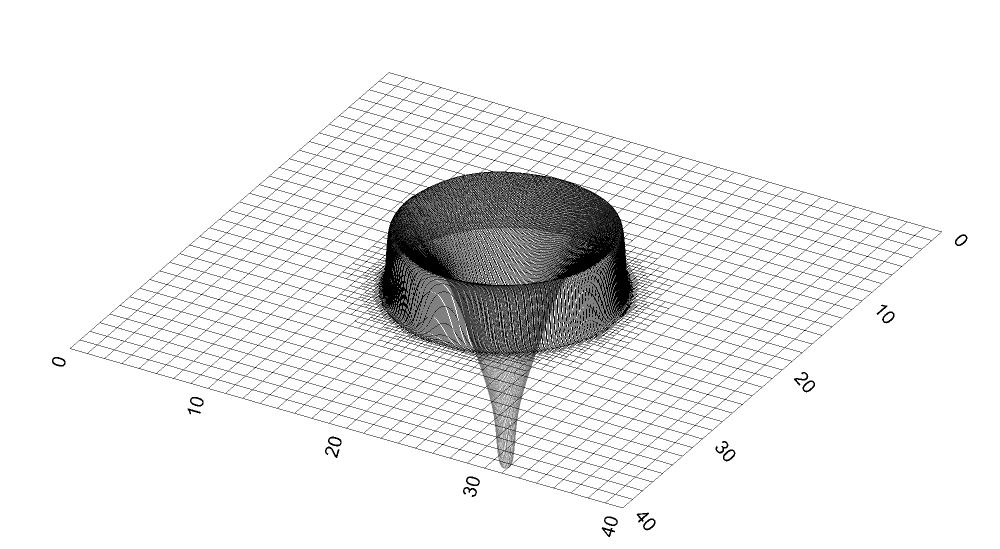}}
\caption{\sf Example 1: Initial and computed water surface $w(x,y,t)$ (left column) and corresponding quadtree grids (right column) for $t=0$, $1$,
$2$, $3$, and $4$ (from top to bottom).\label{fig:3}}
\end{figure}

Furthermore, we take $m=9$ refinement levels of the quadtree grid and set $C_{w,\,\rm seed}= 5 \times 10^{-4}$ in \eqref{3.34}. 
The solution runs until the final time $t=4$. Water surface contours and the respective quadtree grids are illustrated in Figure \ref{fig:3}. 
The quadtree grid starts with 2,134 cells and ends with 35,200 cells at $t=4$. The results in Figure \ref{fig:3} show that the solution follows 
the same evolution in comparison with the ones obtained in \cite{EleuterioFToro2001, Jiang2011}.

\subsection*{Example 2 --- Dam break with density discontinuity over a hump}
This example is based on the benchmark in \cite{Burguete2008, Jiang2011}. We show the capability of the central-upwind 
quadtree scheme to maintain the well-balanced property, symmetry, and to generate adaptive grids. We consider the computational domain 
to be $[0,2] \times [0,1]$ with the following initial conditions:

$$
\rho(x,y,0)=\left\{\begin{array}{lc}997,&~x < 1,\\1200,&~x\geq 1,\end{array}\right. \quad
u(x,y,0)=v(x,y,0)\equiv0, \quad w(x,y,0) \equiv 1,
$$
and the given bottom topography
$$
B(x,y)=0.8e^{-5(x-0.9)^2-50(y-0.5)^2}.
$$

We set $C_{w,\,\rm seed}= 10^{-2}$ and $C_{\rho,\,\rm seed}= 10$ in \eqref{3.34} and \eqref{3.35}, and $m=8$. A solid wall boundary 
condition is used at the top and bottom boundaries. For the sake of simplicity, we set the left and the right boundaries 
to Dirichlet boundary conditions. We run the solution up to the final time $t=0.8$ with the non-well-balanced and 
well-balanced schemes. For the non-well-balanced scheme, the source approximations read as: 

\begin{align*}
&\xbar S_{j,k}^{\,(2)}=
-\frac{g \xbar {(h\rho)}_{j,k}}{ \rho_\circ \, \dx_{j,k}}\left[\frac{B_{\jph,\kph}+B_{\jph,\kmh}}{2}-\frac{B_{\jmh,\kph}+B_{\jmh,\kmh}}{2}\right],\\
&\xbar S_{j,k}^{\,(3)}=
-\frac{g \xbar {(h\rho)}_{j,k}}{ \rho_\circ \, \dx_{j,k}}\left[\frac{B_{\jph,\kph}+B_{\jph,\kmh}}{2}-\frac{B_{\jmh,\kph}+B_{\jmh,\kmh}}{2}\right].
\end{align*}

\begin{figure}[ht!]
\centerline{\includegraphics[height=3.8cm]{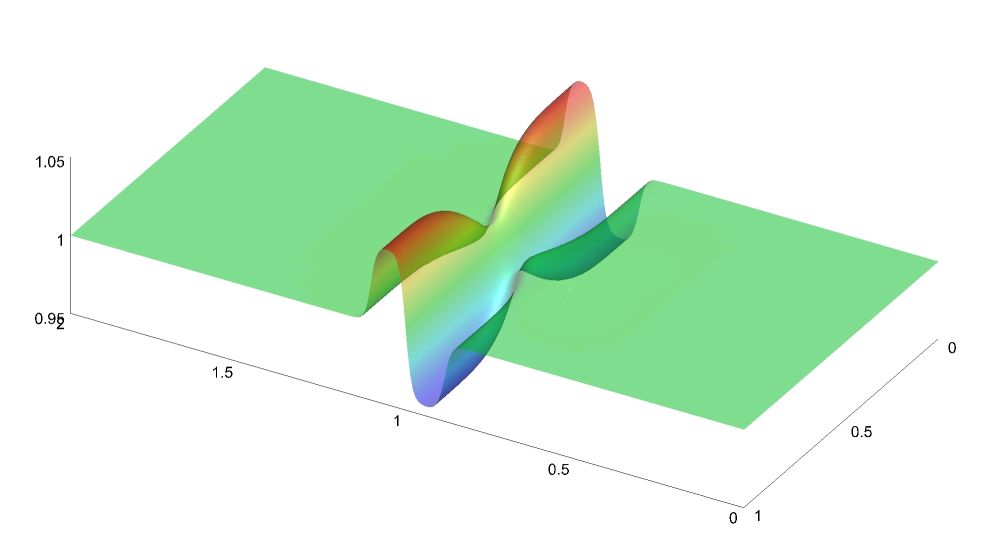}\hspace*{1cm}\includegraphics[height=3.8cm]{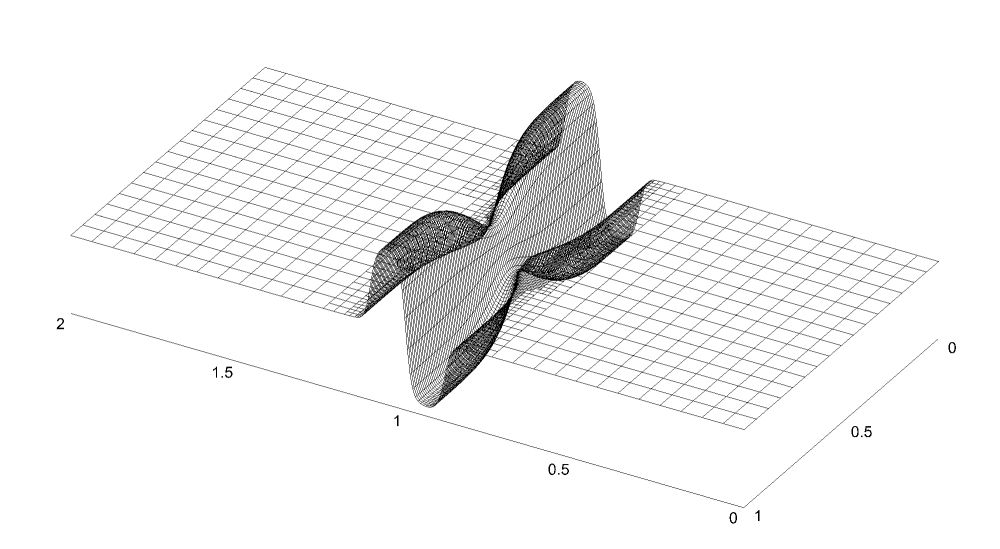}}
\vspace*{0.25cm}
\centerline{\includegraphics[height=3.8cm]{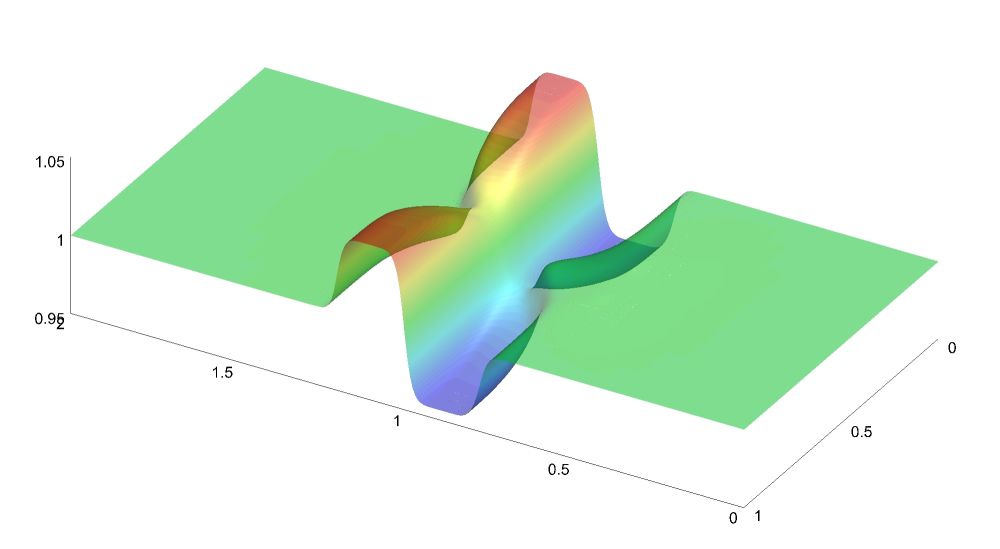}\hspace*{1cm}\includegraphics[height=3.8cm]{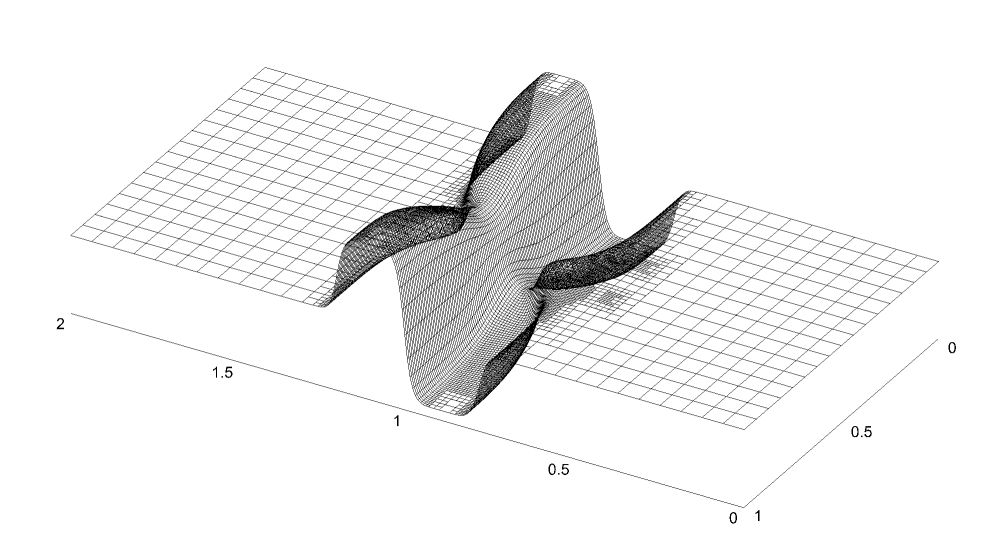}}
\vspace*{0.25cm}
\centerline{\includegraphics[height=3.8cm]{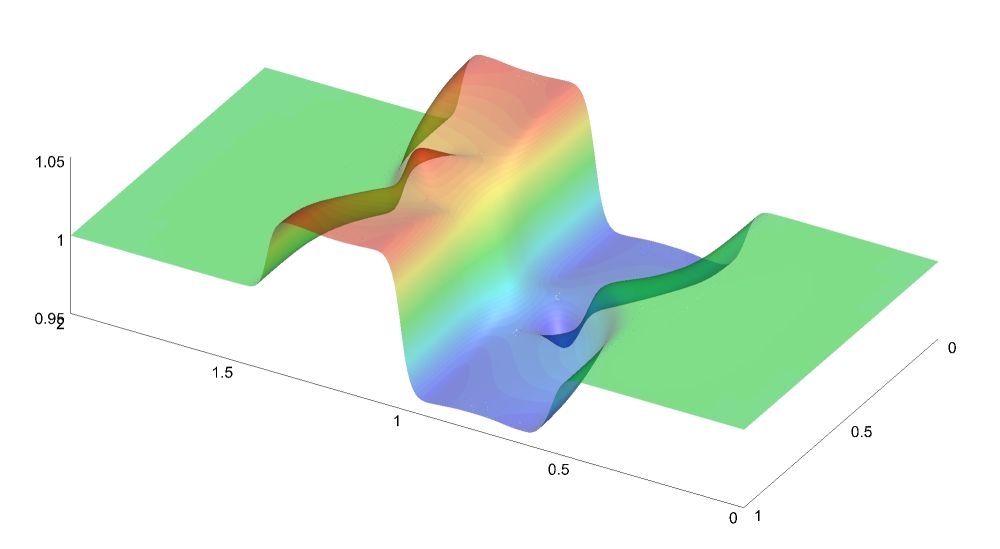}\hspace*{1cm}\includegraphics[height=3.8cm]{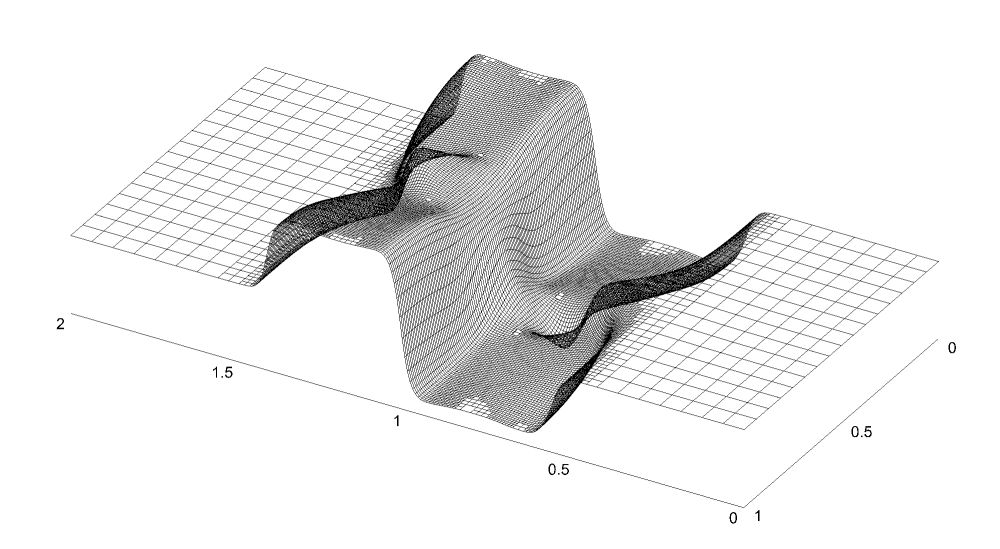}}
\vspace*{0.25cm}
\centerline{\includegraphics[height=3.8cm]{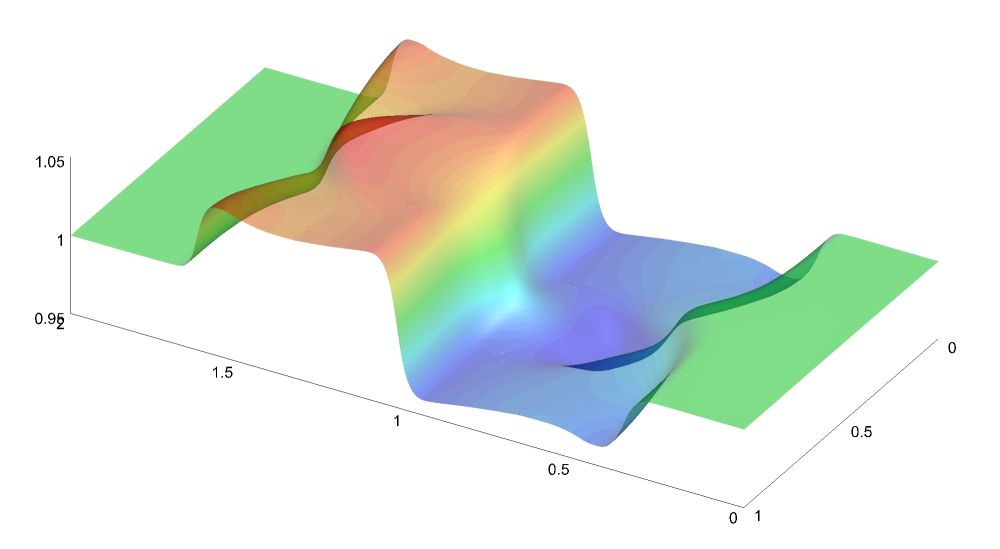}\hspace*{1cm}\includegraphics[height=3.8cm]{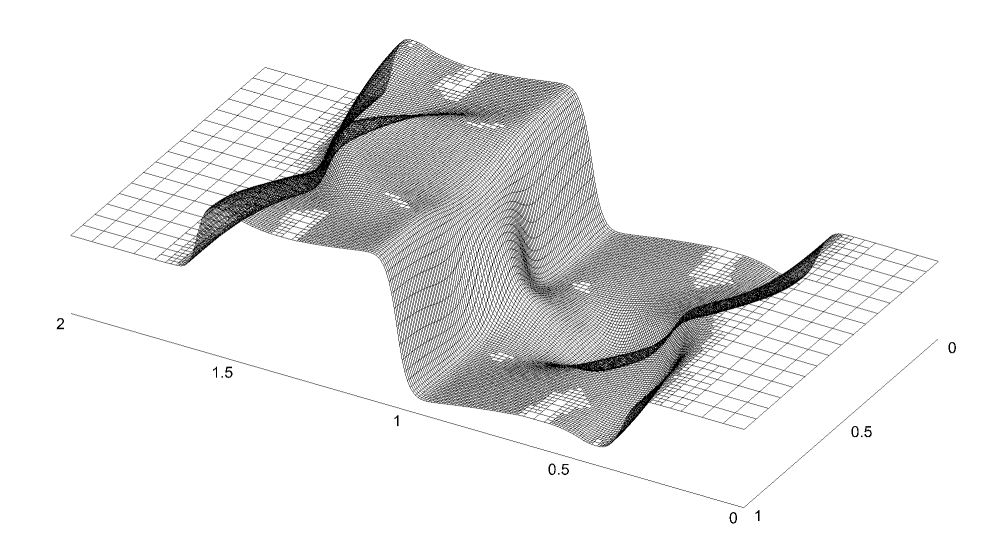}}
\vspace*{0.25cm}
\centerline{\includegraphics[height=3.8cm]{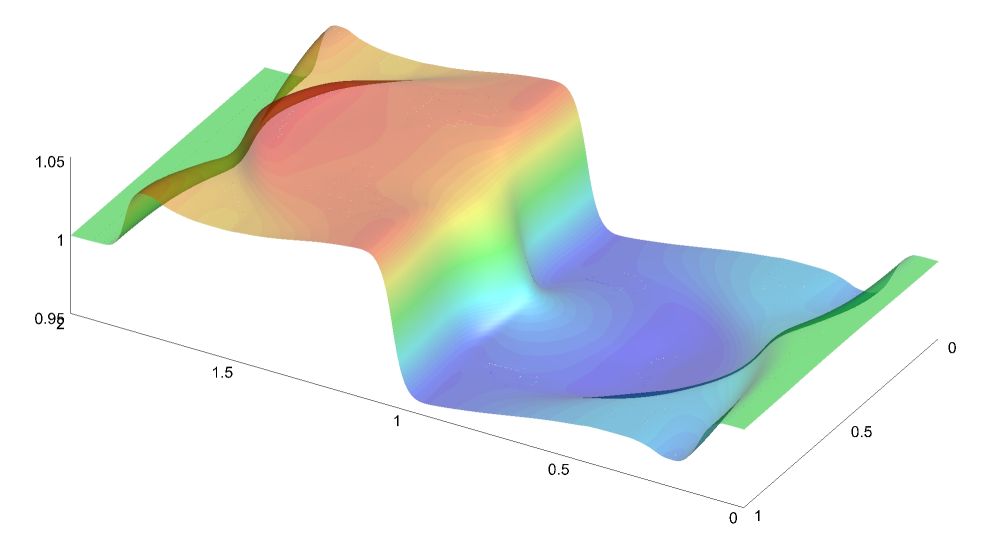}\hspace*{1cm}\includegraphics[height=3.8cm]{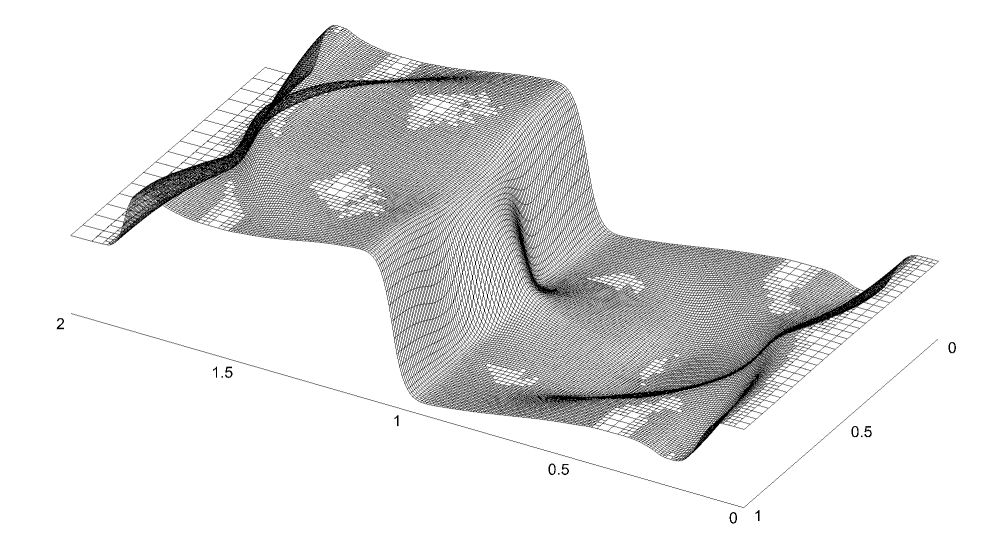}}
\caption{\sf Example 2: Computed water surface $w(x,y,t)$ (left column) and corresponding quadtree grids (right column) for $t=0.1$, $0.2$,
$0.4$, $0.6$, and $0.8$ (from top to bottom) obtained using the non-well-balanced scheme.\label{fig:4}}
\end{figure}

\begin{figure}[ht!]
\centerline{\includegraphics[height=3.8cm]{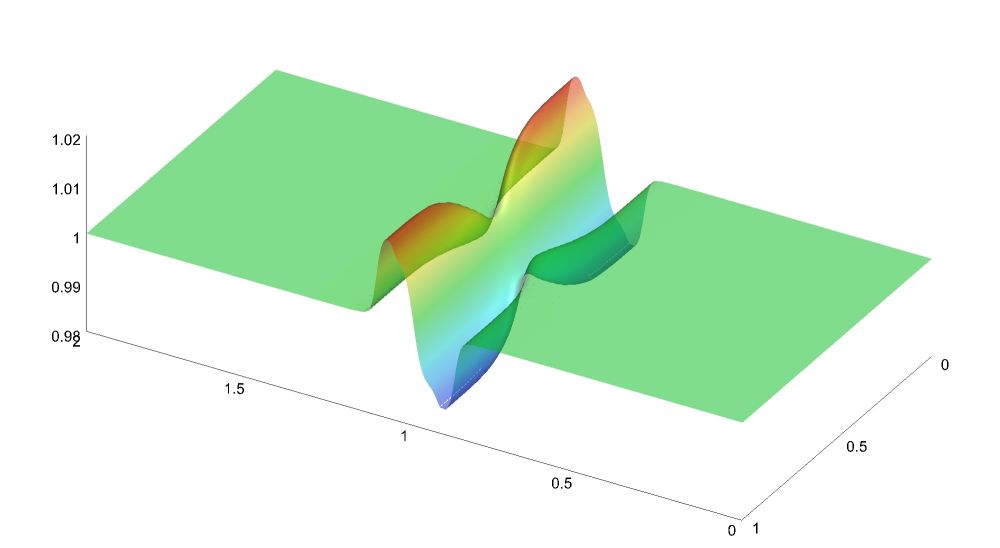}\hspace*{1cm}\includegraphics[height=3.8cm]{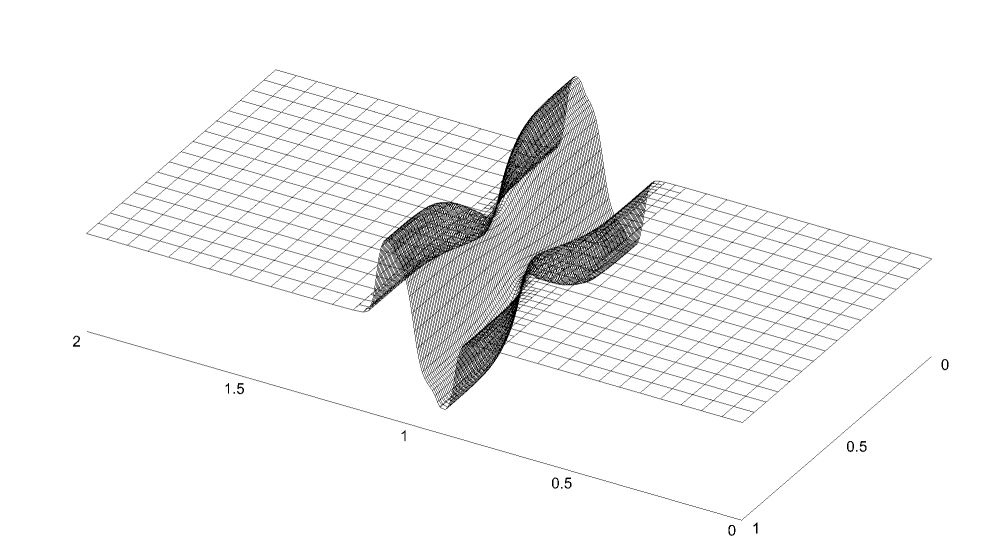}}
\vspace*{0.25cm}
\centerline{\includegraphics[height=3.8cm]{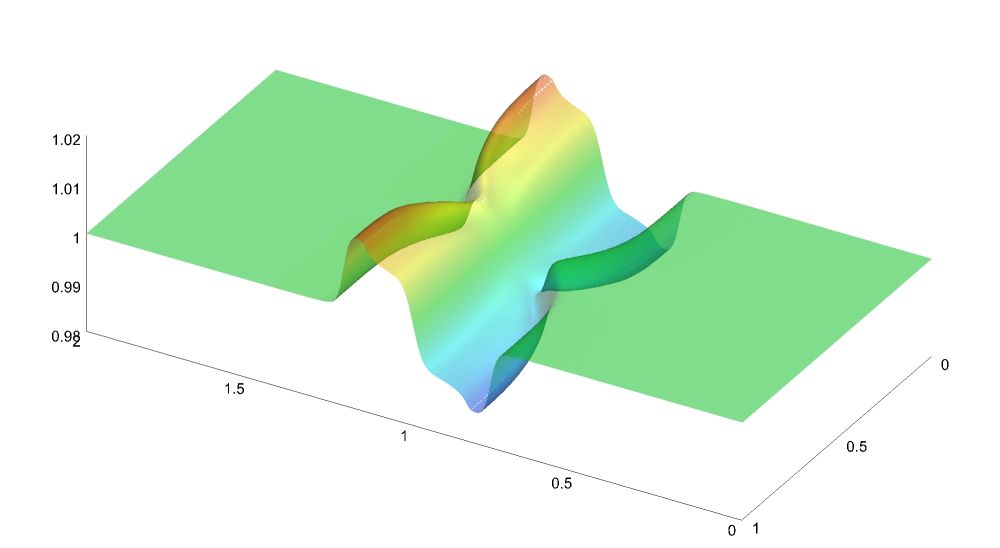}\hspace*{1cm}\includegraphics[height=3.8cm]{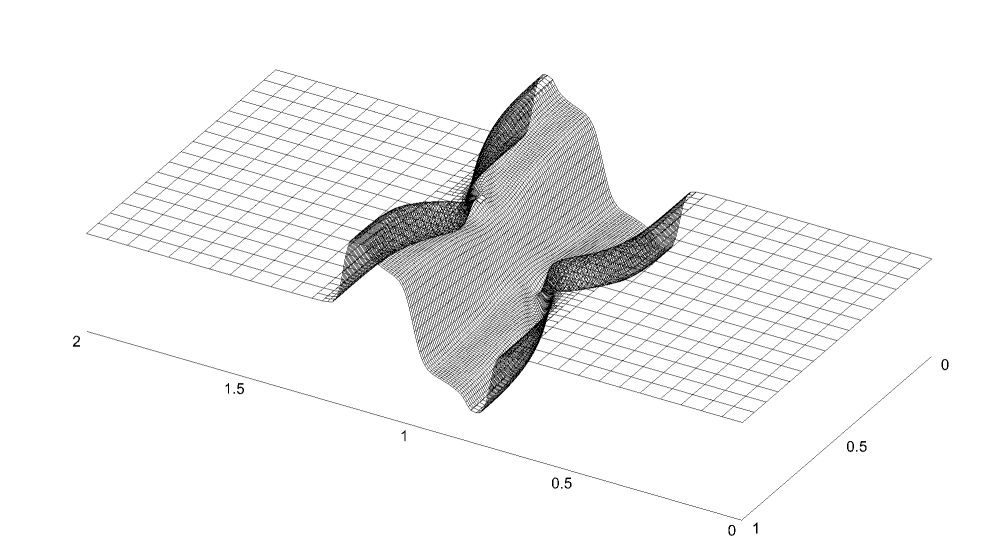}}
\vspace*{0.25cm}
\centerline{\includegraphics[height=3.8cm]{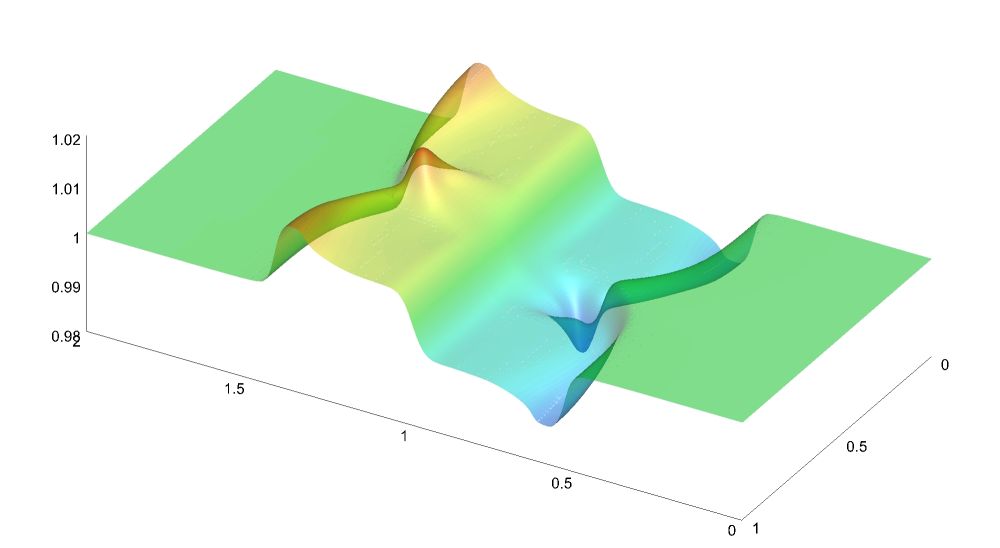}\hspace*{1cm}\includegraphics[height=3.8cm]{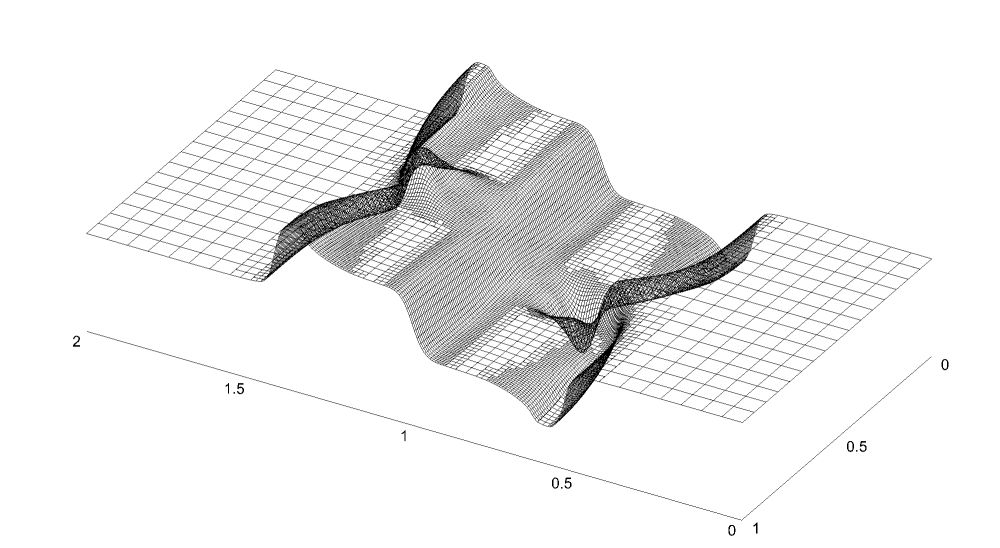}}
\vspace*{0.25cm}
\centerline{\includegraphics[height=3.8cm]{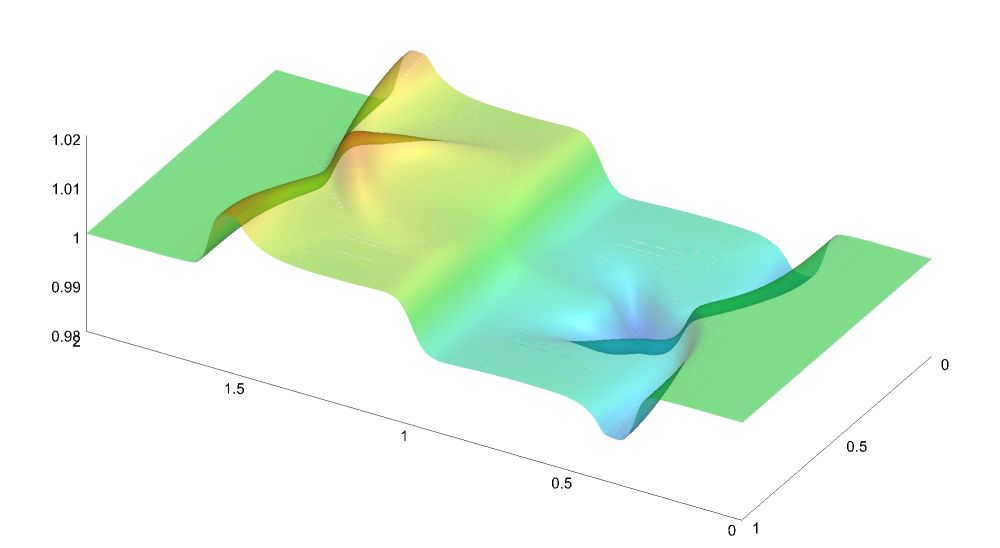}\hspace*{1cm}\includegraphics[height=3.8cm]{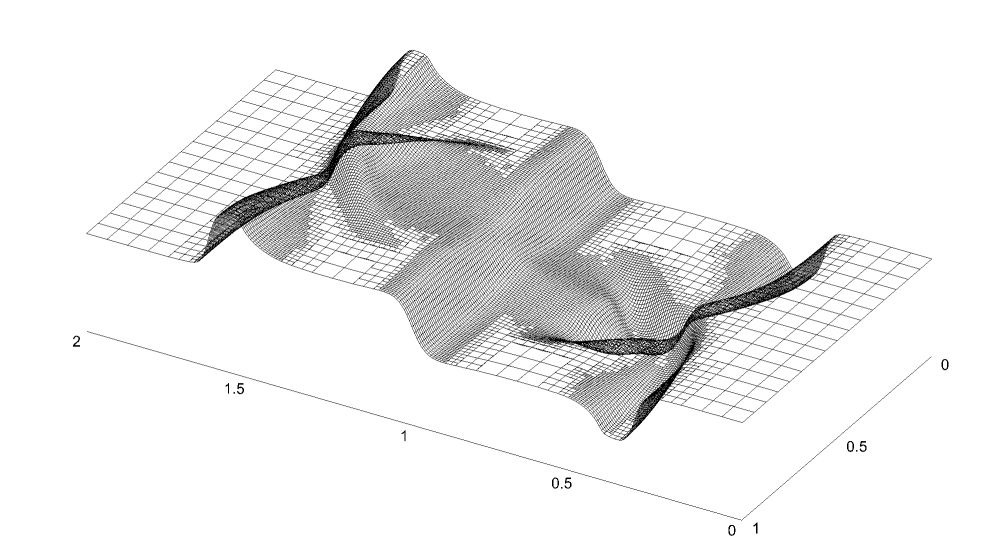}}
\vspace*{0.25cm}
\centerline{\includegraphics[height=3.8cm]{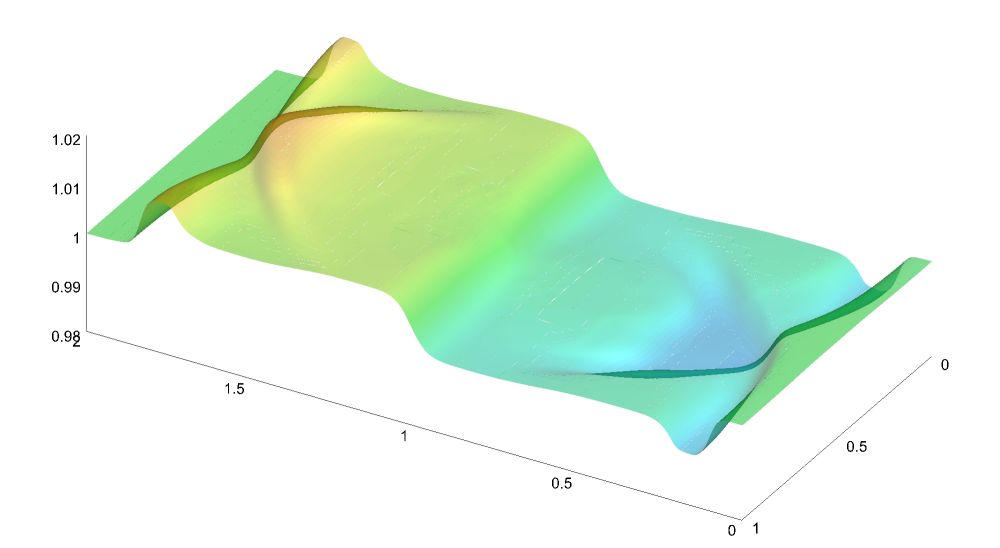}\hspace*{1cm}\includegraphics[height=3.8cm]{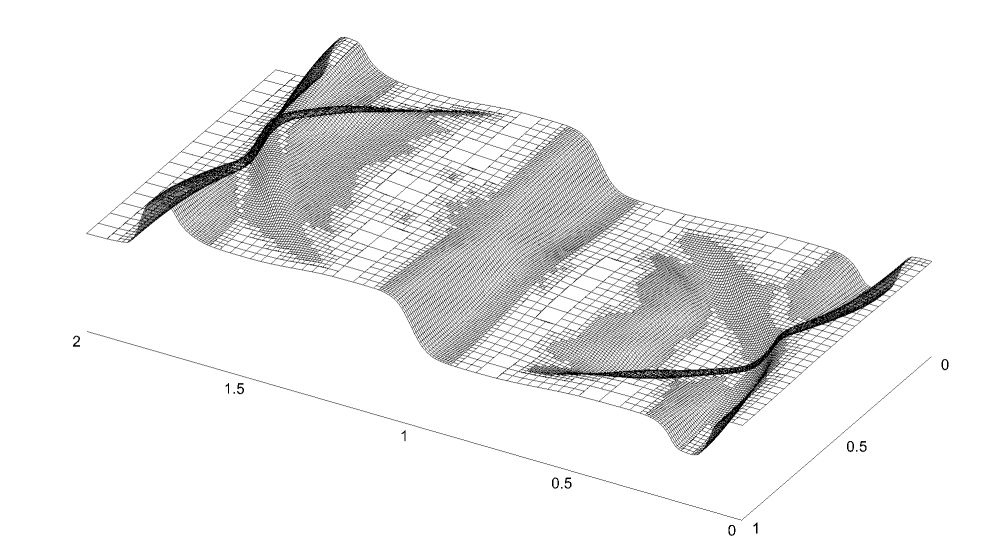}}
\caption{\sf Example 2: Computed water surface $w(x,y,t)$ (left column) and corresponding quadtree grids (right column) for $t=0.1$, $0.2$,
$0.4$, $0.6$, and $0.8$ (from top to bottom) obtained using the well-balanced scheme.\label{fig:5}}
\end{figure}

Water surface contours and the respective quadtree grids of the solution of the non-well-balanced scheme are 
demonstrated in Figure \ref{fig:4}. The quadtree grid starts with 1,184 cells and ends with 27,749 cells at $t=0.8$. 
The well-balanced solution is presented in Figure \ref{fig:5} where the quadtree has a minimum of 1,184 and a maximum of 18,497 cells.
Utilizing the well-balanced scheme reduces the number of cells in the quadtree grid; thereby, the computational 
cost is reduced and unphysical oscillations are eliminated.

\subsection*{Example 3 --- Small perturbations of a stationary steady-state solution}
This numerical example tests the capability of the proposed scheme to capture small perturbations of a steady state solution 
\cite{Bryson2011,Bryson2005, Ghazizadeh2019, Kurganov2002,LeVeque1998,Liu2015}. 
We choose a computational domain $[-2,2]\times[0,1]$ to prevent complicated boundary conditions. The initial conditions are
\begin{align*}
&w(x,y,0)=\left\{\begin{array}{lc}1.01,&~0.05<x<0.15,\\1,&~\mbox{otherwise},\end{array}\right.\qquad u(x,y,0)=v(x,y,0)\equiv0,\\
&\rho(x,y,0)=\left\{\begin{array}{lc}1007,&~0.05<x<0.15,\\\rho_\circ,&~\mbox{otherwise},\end{array}\right.
\end{align*}
and the given bottom topography function in Example 2.

\begin{figure}[ht!]
\centerline{\includegraphics[height=3.6cm]{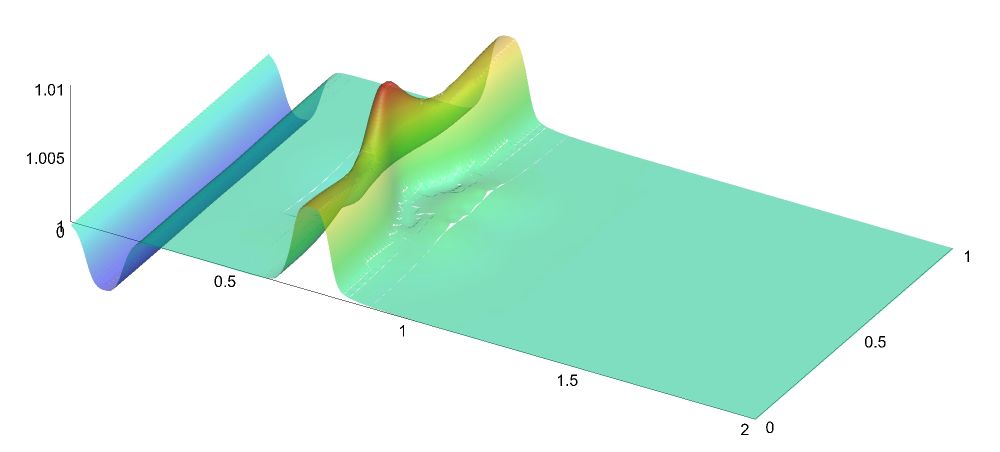}\hspace*{0.5cm}\includegraphics[height=3.6cm]{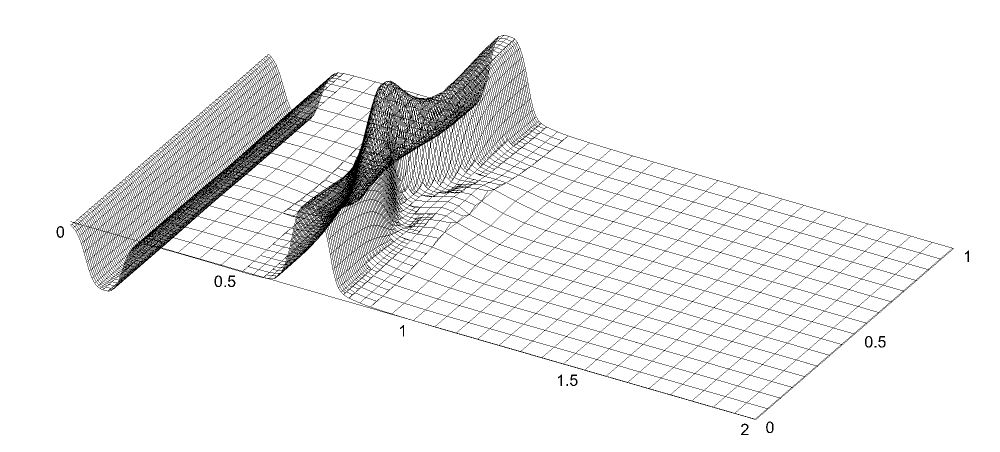}}
\vspace*{0.25cm}
\centerline{\includegraphics[height=3.6cm]{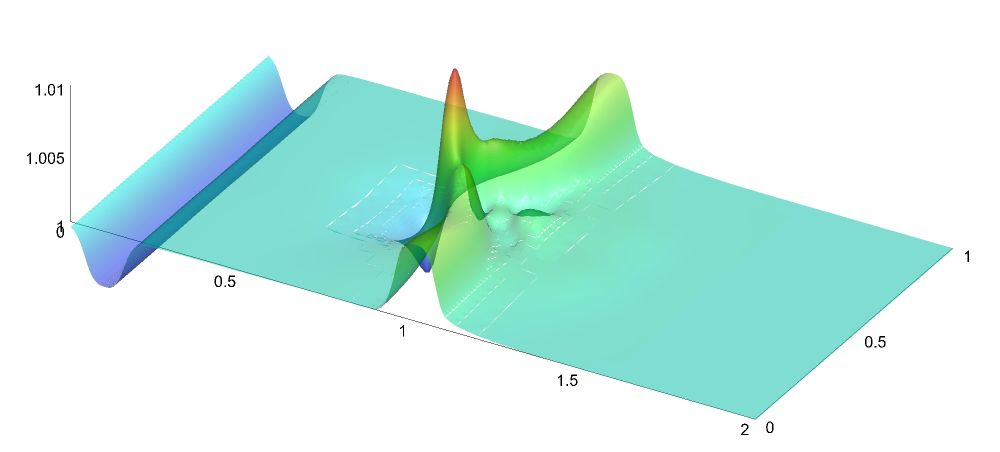}\hspace*{0.5cm}\includegraphics[height=3.6cm]{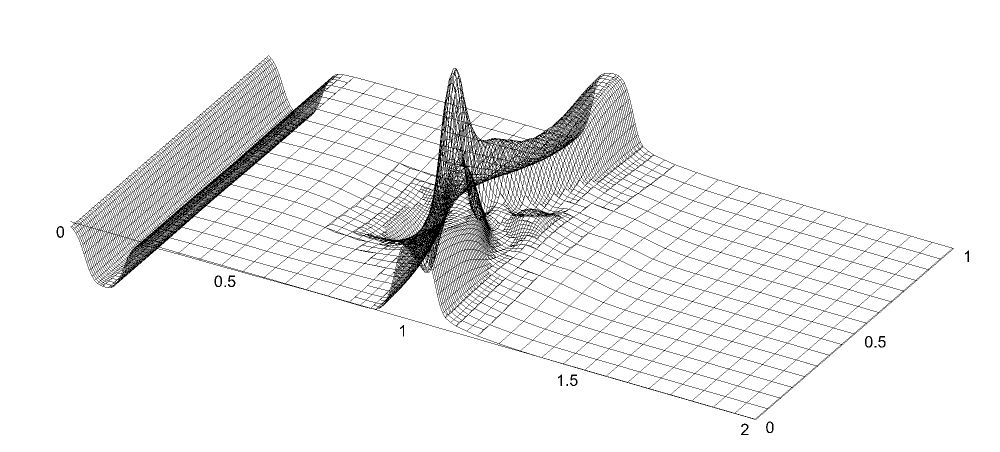}}
\vspace*{0.25cm}
\centerline{\includegraphics[height=3.6cm]{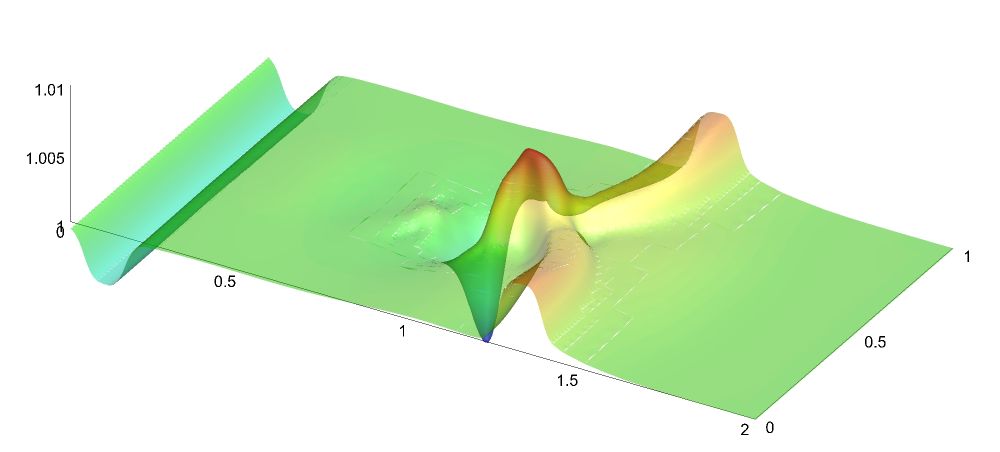}\hspace*{0.5cm}\includegraphics[height=3.6cm]{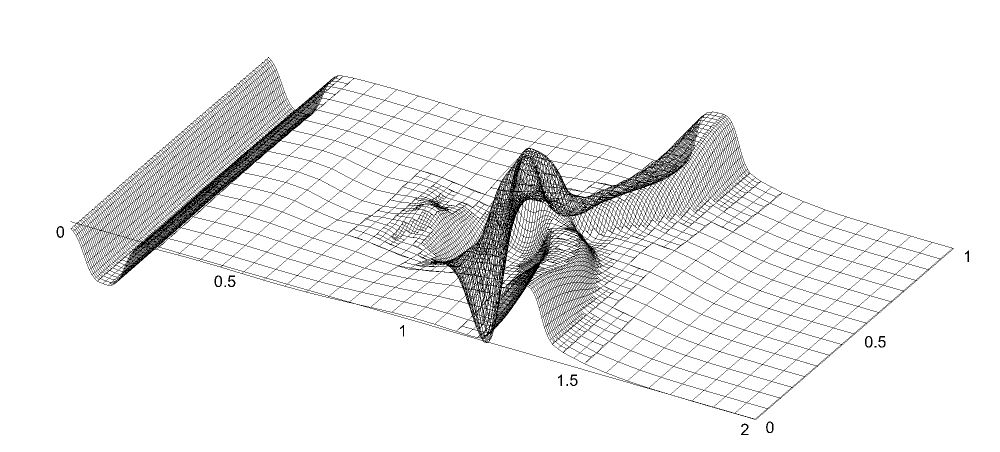}}
\vspace*{0.25cm}
\centerline{\includegraphics[height=3.6cm]{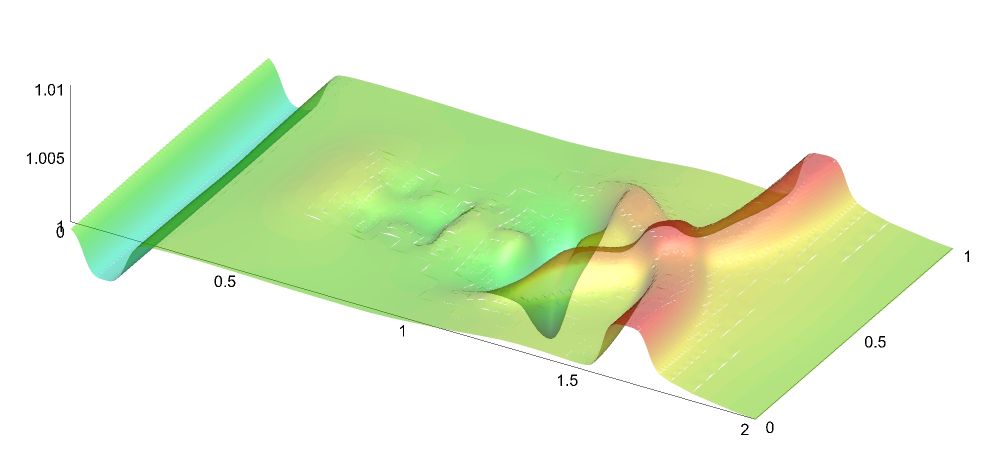}\hspace*{0.5cm}\includegraphics[height=3.6cm]{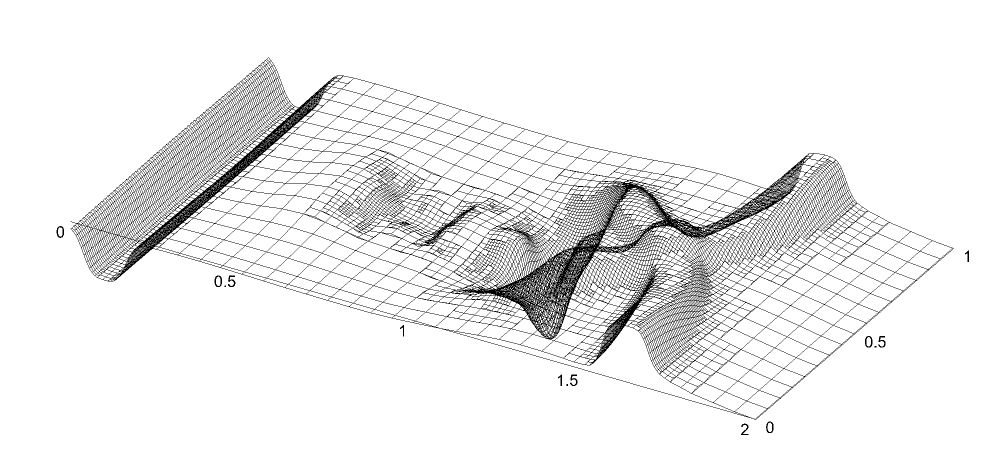}}
\vspace*{0.25cm}
\centerline{\includegraphics[height=3.6cm]{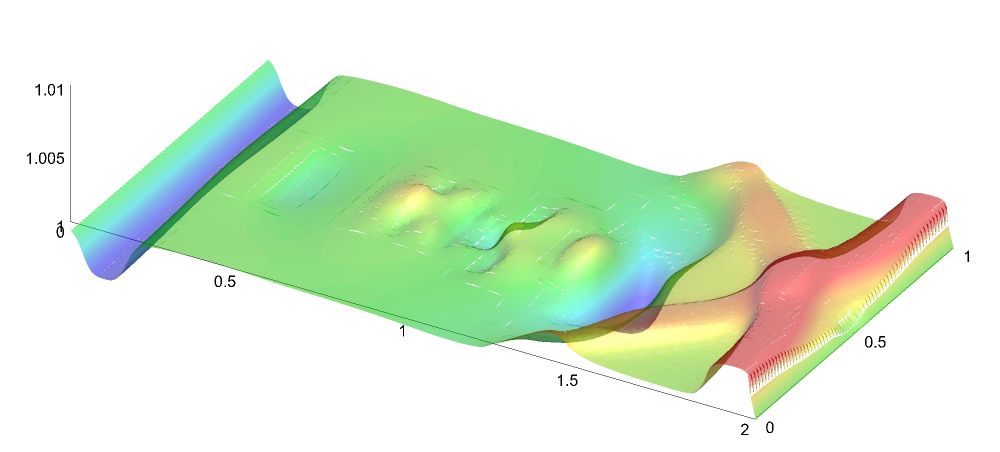}\hspace*{0.5cm}\includegraphics[height=3.6cm]{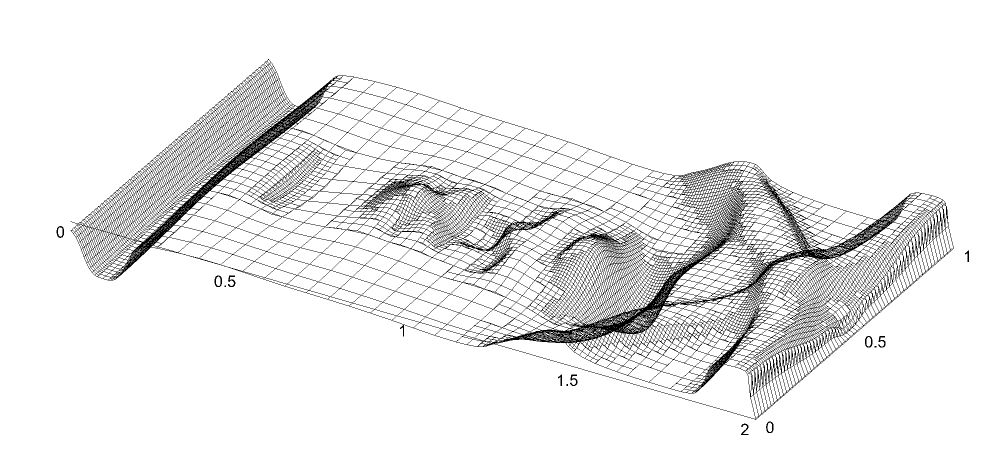}}
\caption{\sf Example 3: Computed water surface $w(x,y,t)$ (left column) and corresponding quadtree grids (right column) for $t=0.6$, $0.9$,
$1.2$, $1.5$, and $1.8$ (from top to bottom) obtained using the non-well-balanced scheme.\label{fig:6}}
\end{figure}

\begin{figure}[ht!]
\centerline{\includegraphics[height=3.6cm]{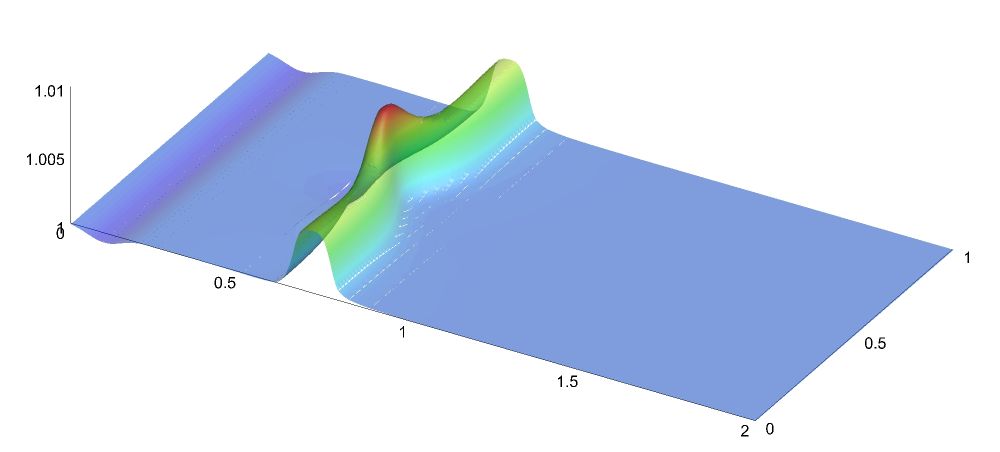}\hspace*{0.5cm}\includegraphics[height=3.6cm]{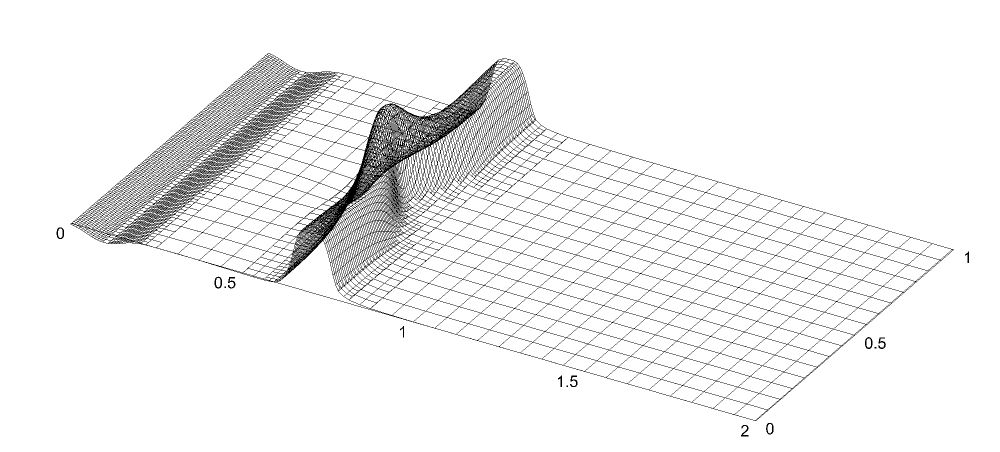}}
\vspace*{0.25cm}
\centerline{\includegraphics[height=3.6cm]{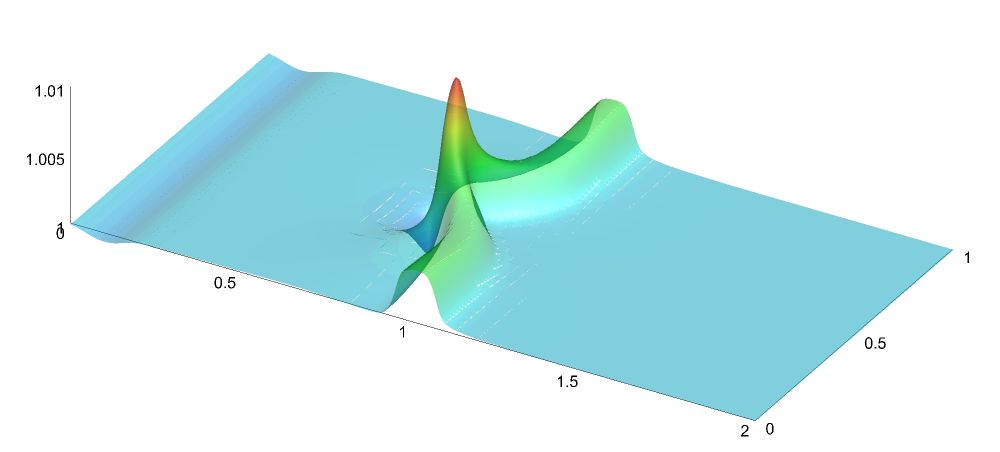}\hspace*{0.5cm}\includegraphics[height=3.6cm]{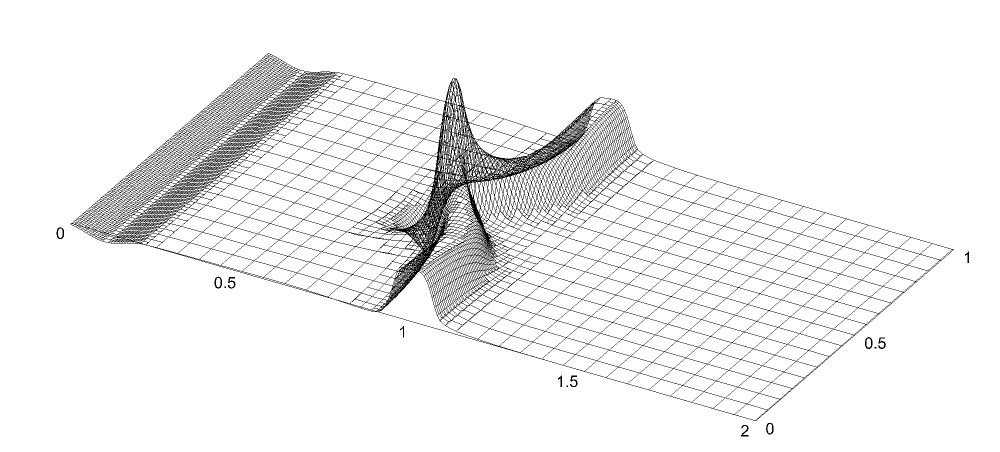}}
\vspace*{0.25cm}
\centerline{\includegraphics[height=3.6cm]{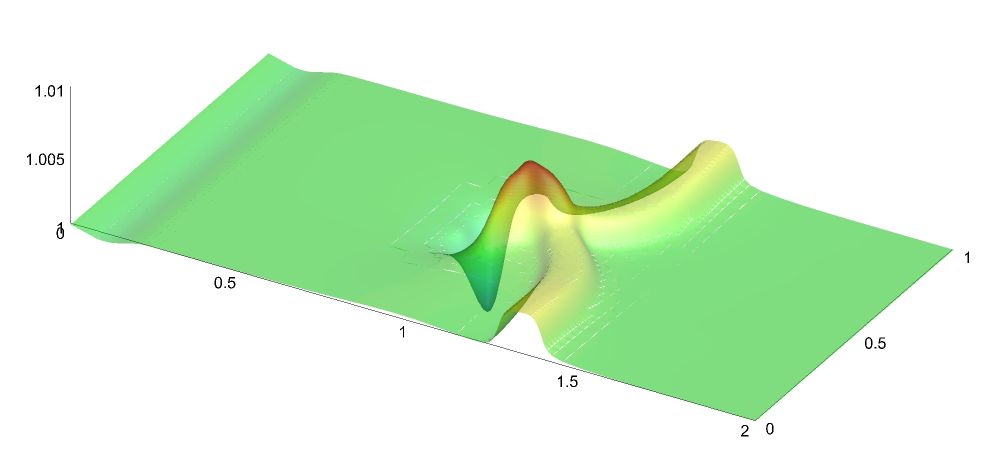}\hspace*{0.5cm}\includegraphics[height=3.6cm]{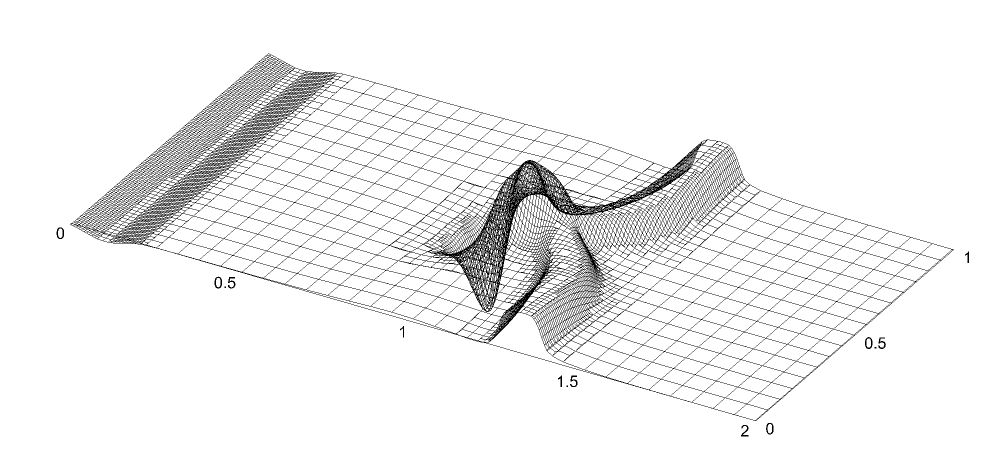}}
\vspace*{0.25cm}
\centerline{\includegraphics[height=3.6cm]{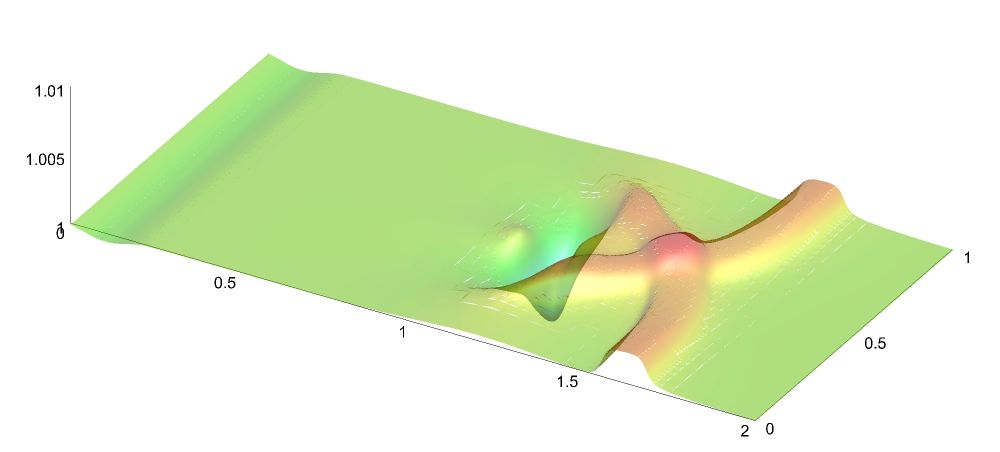}\hspace*{0.5cm}\includegraphics[height=3.6cm]{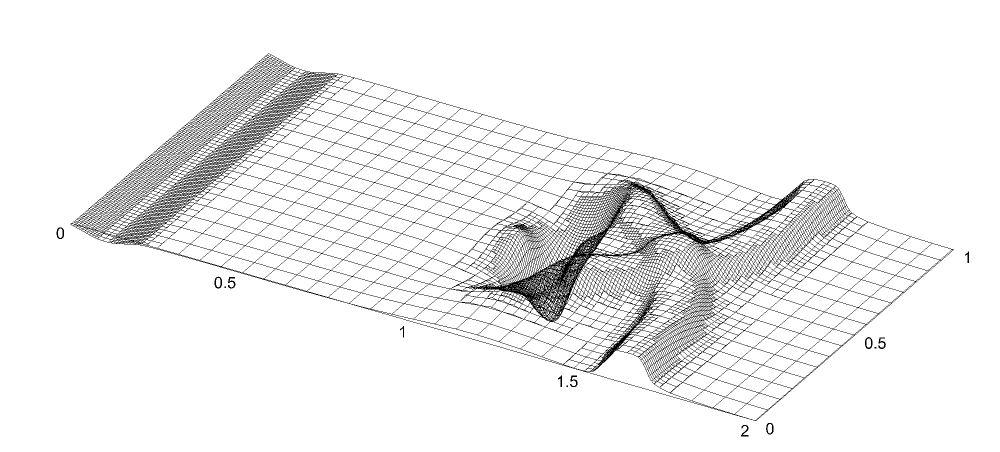}}
\vspace*{0.25cm}
\centerline{\includegraphics[height=3.6cm]{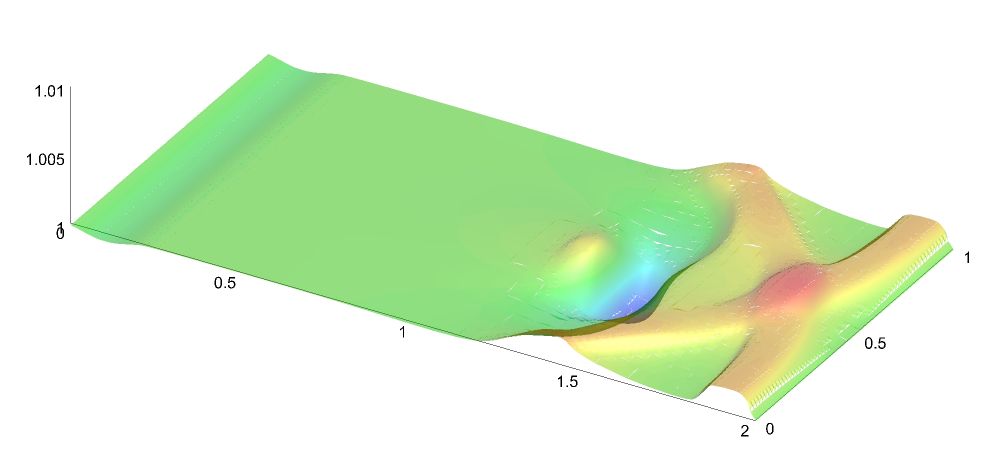}\hspace*{0.5cm}\includegraphics[height=3.6cm]{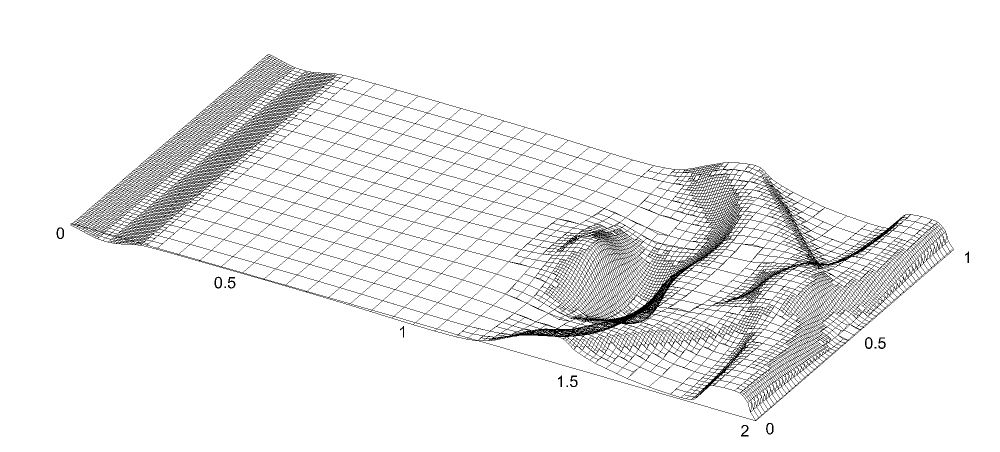}}
\caption{\sf Example 3: Computed water surface $w(x,y,t)$ (left column) and corresponding quadtree grids (right column) for $t=0.6$, $0.9$,
$1.2$, $1.5$, and $1.8$ (from top to bottom) obtained using the well-balanced scheme.\label{fig:7}}
\end{figure}

We set boundary conditions similar to Example 2 for this test. $m=9$ is taken, and we set $C_{w, \,\rm seed}=10^{-2}$ and $C_{\rho, \,\rm seed}=10$ 
in \eqref{3.34} and \eqref{3.35}. The non-well-balanced solution is computed until the final time $t=1.8$ and plot 
the snapshots of $w$ (left) and the quadtree grids (right) at times $t=0.6$, $0.9$, $1.2$, $1.5$, and $1.8$ in the domain of $[0,2]\times[0,1]$ in
Figure \ref{fig:6}. The quadtree grid starts with 2,530 cells and reaches a maximum number of 13,804 cells during the 
time evolution. The well-balanced solution is illustrated in Figure \ref{fig:7}, respectively. In this solution, the number of cells reaches the maximum of 10,222.
Figure \ref{fig:7} demonstrates that the proposed well-balanced central-upwind quadtree scheme accurately captures small
perturbations of the ``lake-at-rest'' steady state and that the symmetry of the solution is preserved.

\subsection{Example 4 --- Sudden contraction with variable density inflow}
The last example is a modification of the example in \cite{Ghazizadeh2019, Hubbard2001}. The purpose of this
example is to demonstrate the positivity-preserving property of the proposed scheme.

We consider an open channel with a sudden contraction. The geometry of the channel is established on its contraction, where
$$
y_b(x)=\left\{\begin{array}{lc}0.5,&x\le1,\\0.4,&\mbox{otherwise}.\end{array}\right.
$$
The computational domain is $[0,3]\times[0.5-y_b(x),{0.5+y}_b(x)]$. Solid wall boundary conditions are imposed at all of the boundaries
except for part of the left inflow boundary, with $u(0,y_i,t)\equiv2$ and $\rho(0,y_i,t)\equiv1007$, where $y_i\in[0.4, 0.6]$. In addition, we set the right boundary to a zero-order extrapolation. The following initial conditions are prescribed:
$$
w(x,y,0)\equiv1,\quad u(x,y,0) = v(x,y,0)\equiv0,\quad \rho(x,y,0)\equiv\rho_\circ. 
$$
In this example, we take $m=8$ refinement levels of the quadtree grid and set $C_{w, \,\rm seed}=2$ and $C_{\rho, \,\rm seed}=20$ in 
\eqref{3.34} and \eqref{3.35}. We compute the solution with the bottom topography given in Example 2, where the water depth at 
the top of the humps is quite shallow, which makes it a good example to test the positivity-preserving property.

\begin{figure}[ht!]
  \centerline{\includegraphics[height=2.8cm]{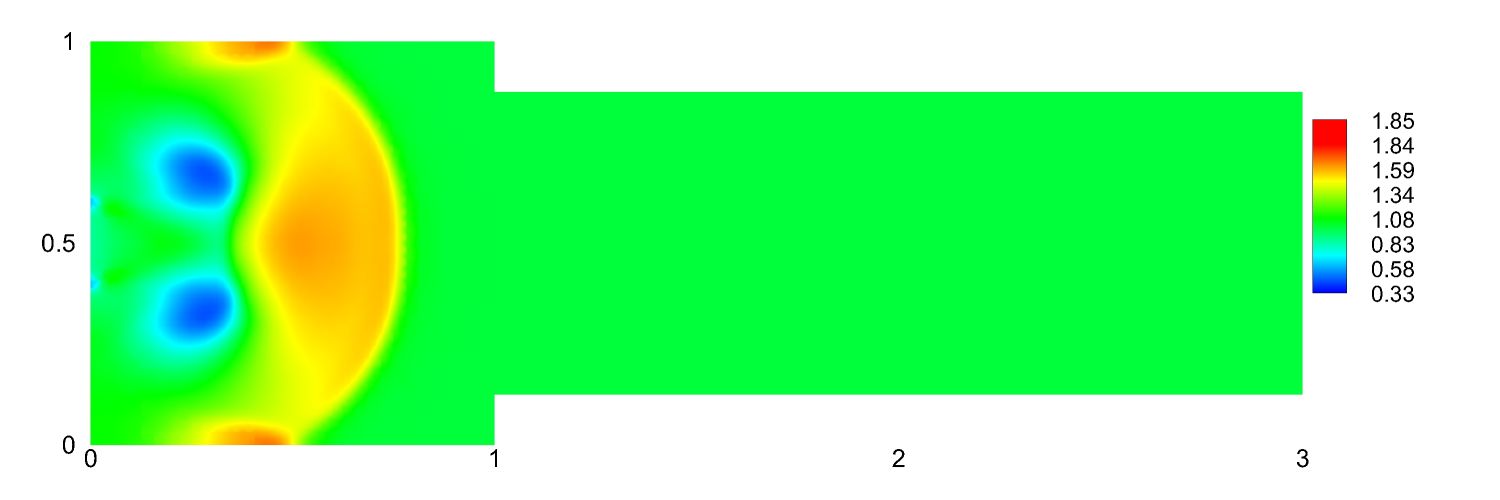}\hspace*{0.1cm}\includegraphics[height=2.8cm]{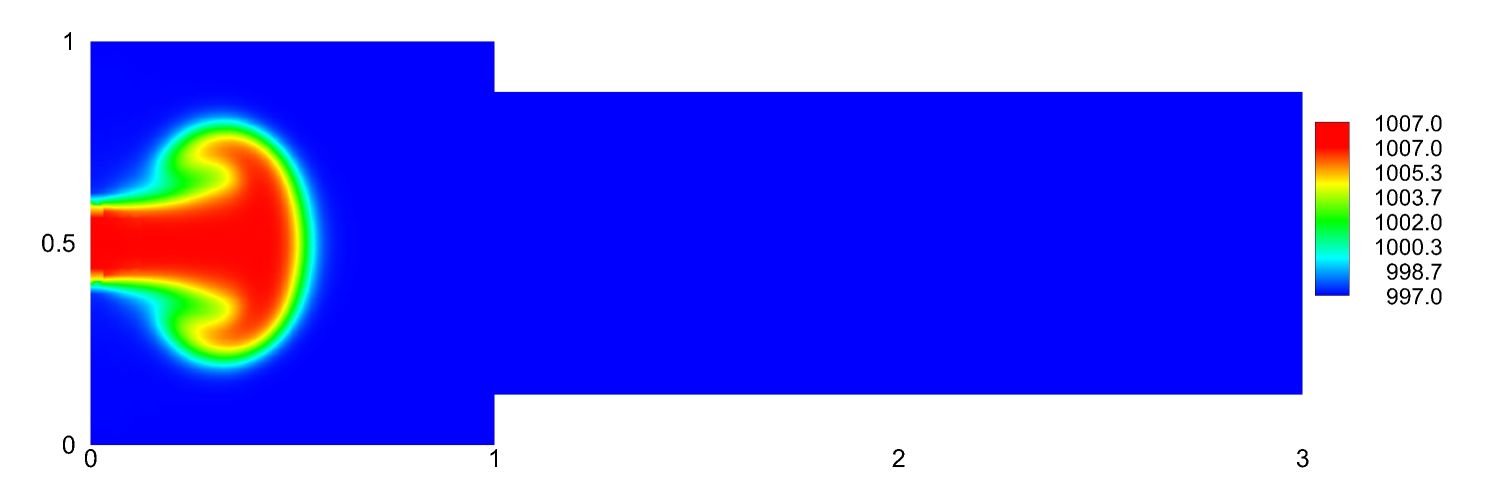}}
  \vspace*{0.25cm}
  \centerline{\includegraphics[height=2.8cm]{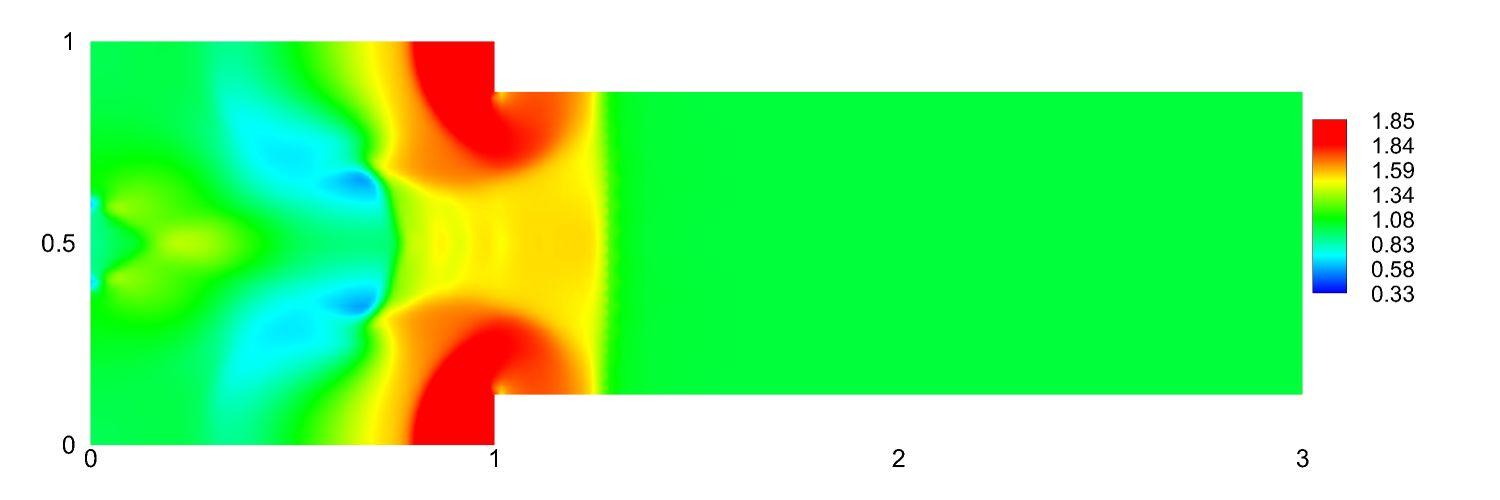}\hspace*{0.1cm}\includegraphics[height=2.8cm]{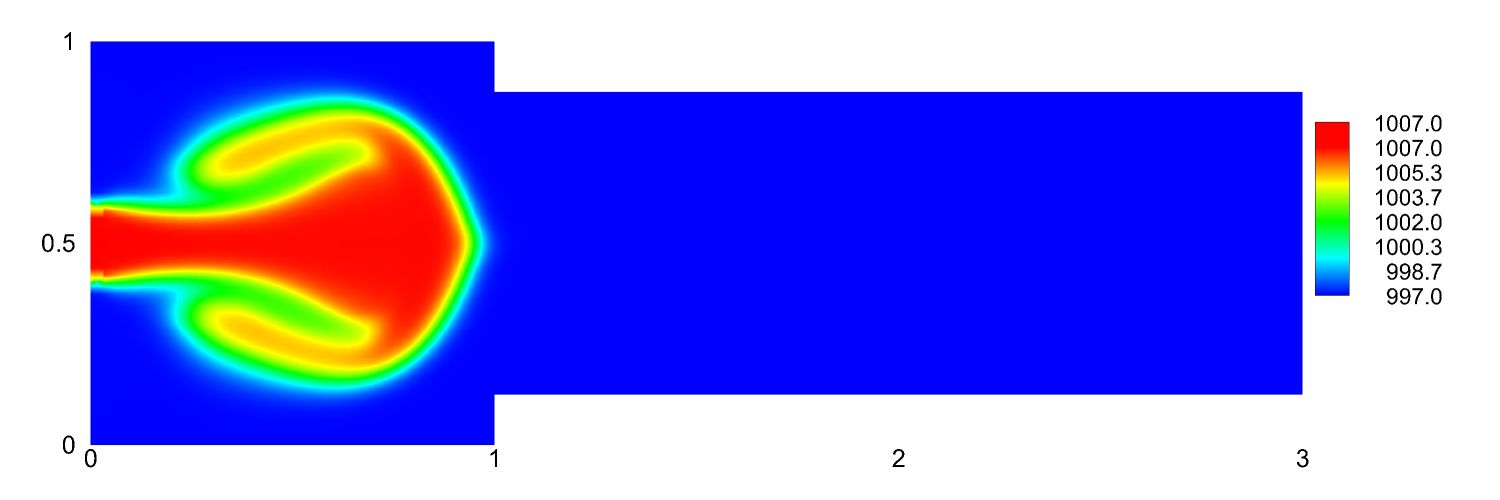}}
  \vspace*{0.25cm}
  \centerline{\includegraphics[height=2.8cm]{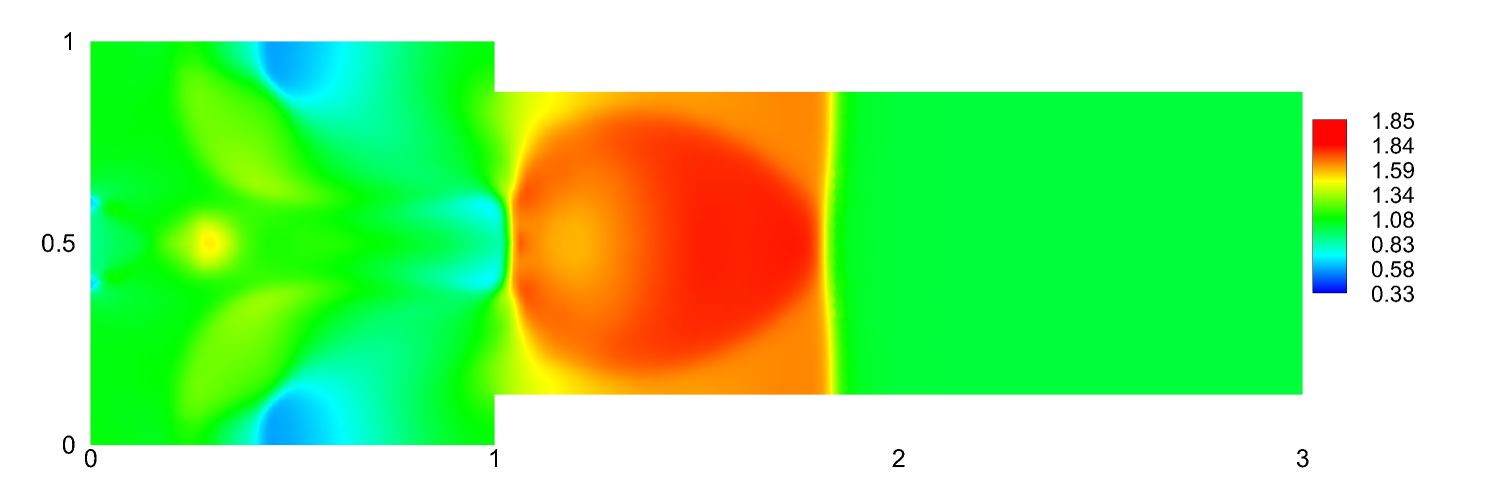}\hspace*{0.1cm}\includegraphics[height=2.8cm]{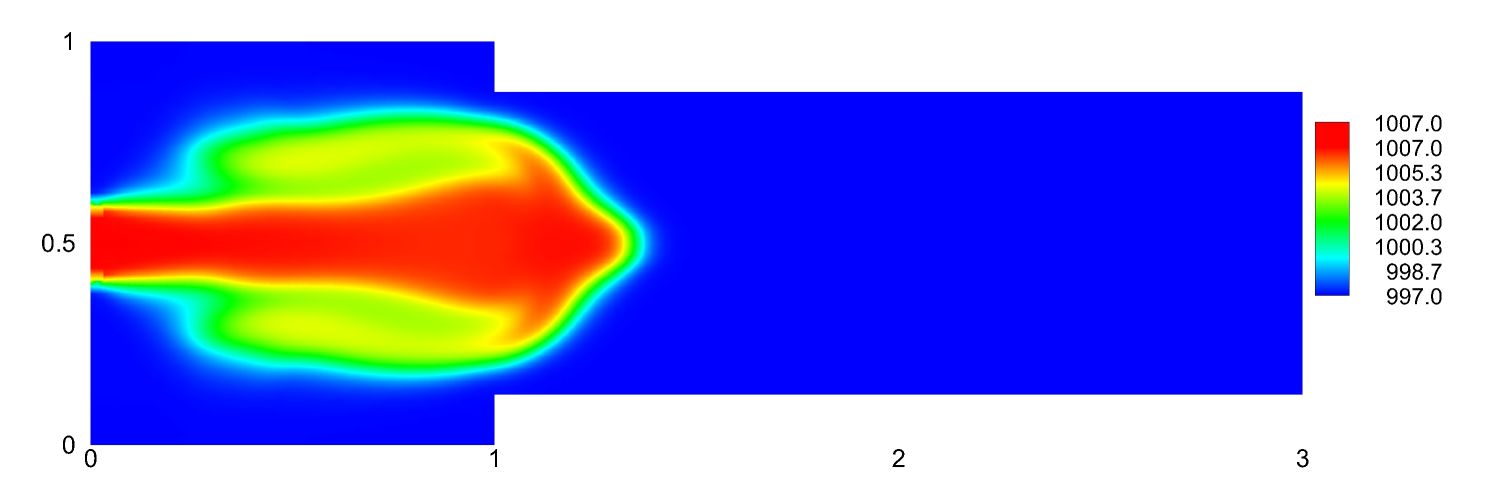}}
  \vspace*{0.25cm}
  \centerline{\includegraphics[height=2.8cm]{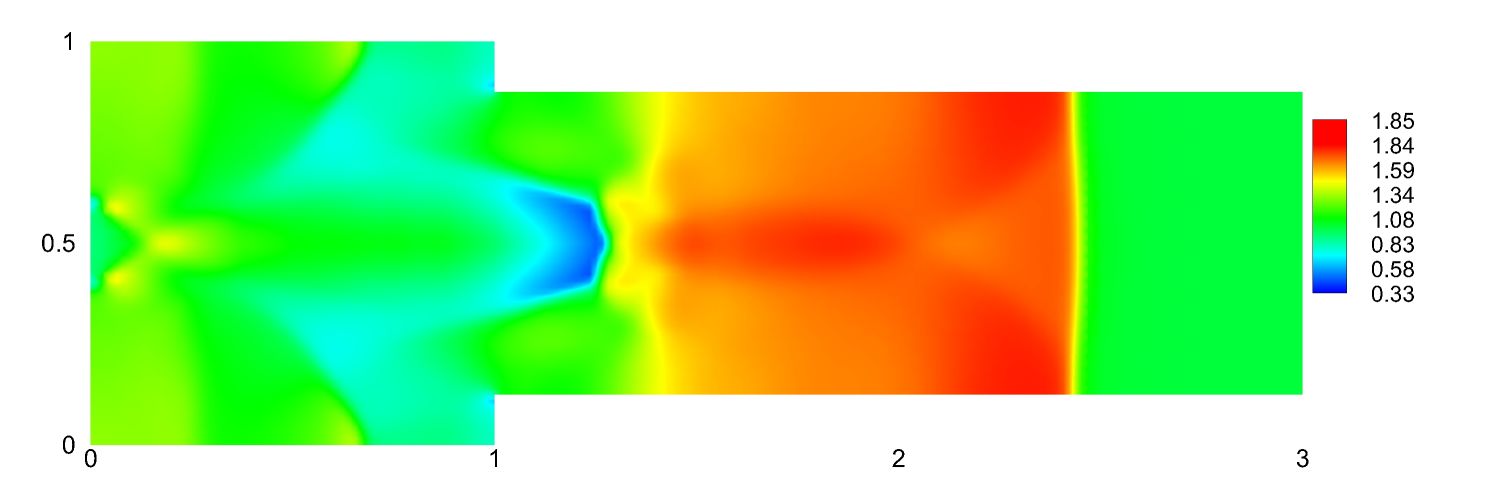}\hspace*{0.1cm}\includegraphics[height=2.8cm]{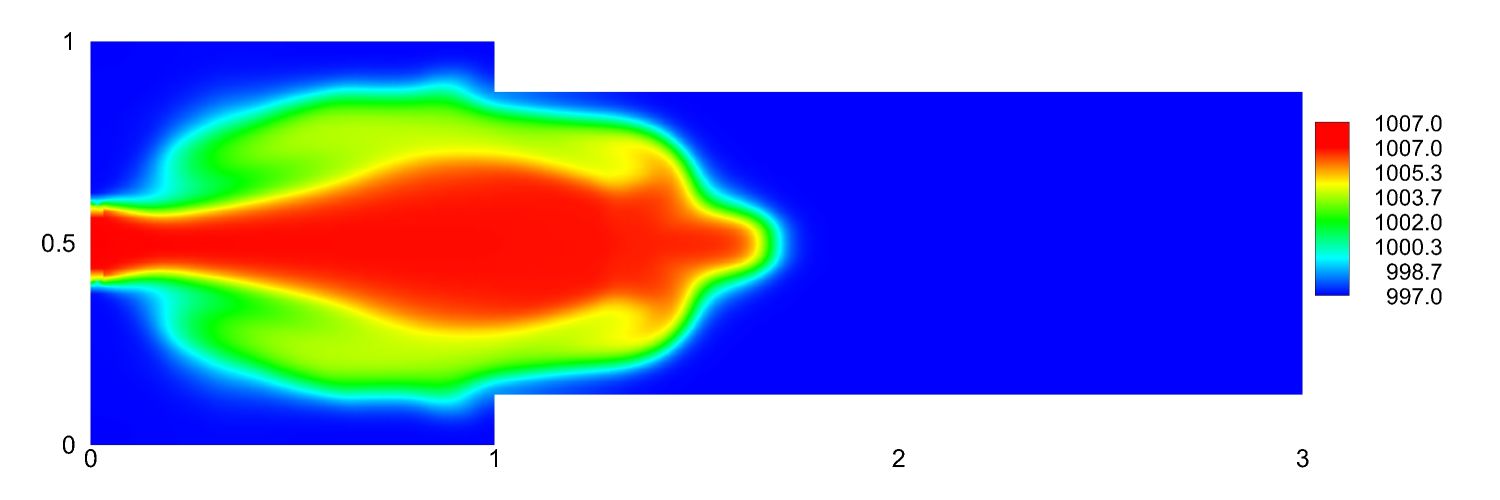}}
  \vspace*{0.25cm}
  \centerline{\includegraphics[height=2.8cm]{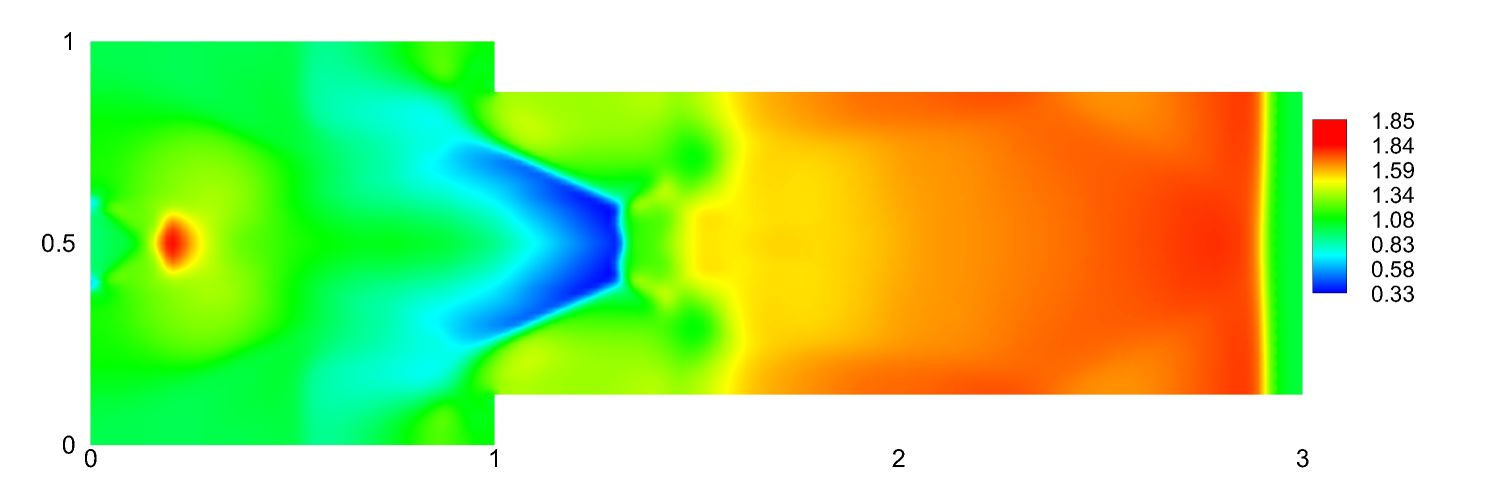}\hspace*{0.1cm}\includegraphics[height=2.8cm]{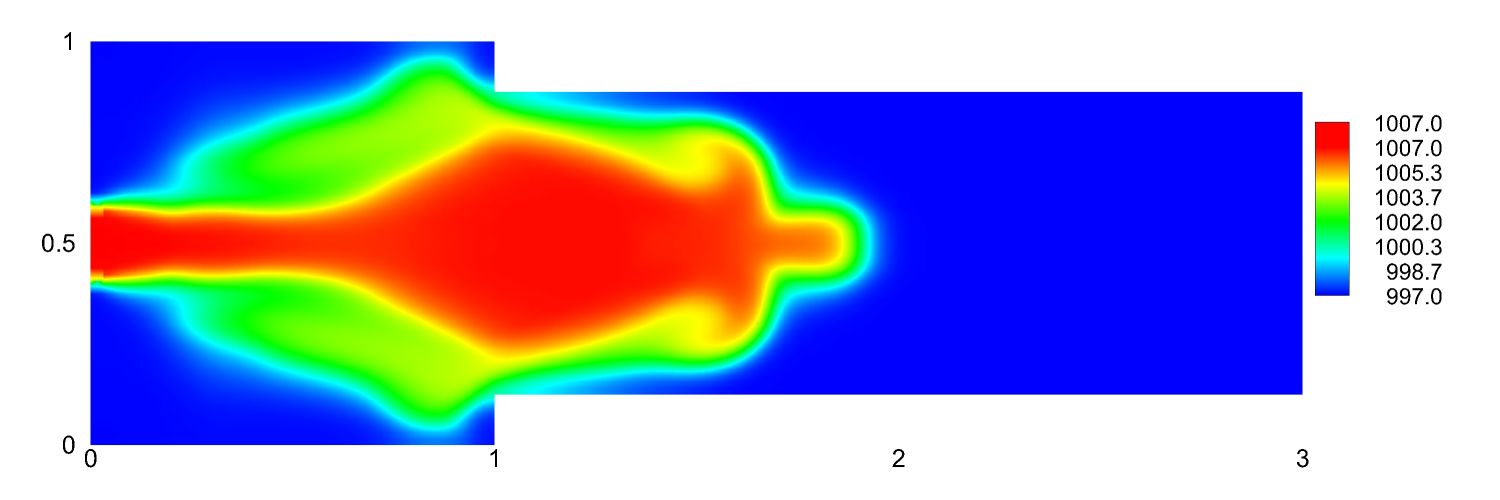}}
  \caption{\sf Example 4: Evolution of water surface $w(x,y,t)$ (left column) and density (right column) for
  $t=0.4$, $0.8$, $1.2$, $1.6$, and $1.9$ (from top to bottom) obtained using the well-balanced scheme.\label{fig:8}}
  \end{figure}

We compute the solution until the final time $t=1.9$ and plot the evolution of $w$ and $\rho$ at times
$t=0.4$, $0.8$, $1.2$, $1.6$, and $1.9$ in Figure \ref{fig:8}. The quadtree grid in this solution starts with a minimum of $298$ and 
reaches a maximum of 9,928 cells. As one can see, the proposed central-upwind quadtree scheme preserves the positivity of 
the computed water depth and density.

\section{Conclusion}\label{S:5}
An adaptive, well-balanced, positivity-preserving central-upwind scheme over quadtree grids for variable density shallow water equations has been presented. Four numerical examples have been used in order to verify the robustness 
and accuracy of the proposed scheme. These tests show symmetry preserving, well-balanced property, positivity-preserving, 
as well as adaptability of the coupled system. The results show that the proposed central-upwind quadtree 
scheme can improve the performance and efficiency of calculations compared to regular Cartesian grids.

\section*{Acknowledgments}
The work of the authors was supported by NSERC grant 210717. The authors warmly thank Carlos Parés from the University of Málaga for providing resources on the well-balanced property 
of the shallow water equations. The authors also thank Yangyang Cao and Philippe LeFloch from the Laboratoire Jacques-Louis Lions 
of Sorbonne  Universit\'e for their useful discussions on the well-balanced
property on conservation laws.




\bibliographystyle{elsarticle-num-names-alphsort}
\bibliography{References}

\end{document}